\documentclass[12pt]{article}
\usepackage{amssymb}
\usepackage{latexsym}
\usepackage{amsfonts}
\usepackage{amsmath}
\newtheorem{theorem}{Theorem}[section]
\newtheorem{lemma}[theorem]{Lemma}
\newtheorem{proposition}[theorem]{Proposition}
\newtheorem{definition}[theorem]{Definition}

\newtheorem{remark}[theorem]{Remark}
\newtheorem{assumption}{Assumption}

\numberwithin{equation}{section}

\usepackage[dvips]{color}

\def\P{\mathbb{P} }
\def\R{\mathbb{R} }

\def\la{\langle}
\def\ra{\rangle}

\begin{document}
\allowdisplaybreaks

\title{\Large\bf
Supercritical Superprocesses:
Proper Normalization and Non-degenerate Strong Limit
}
\author{ \bf  Yan-Xia Ren\footnote{The research of this author is supported by NSFC (Grant No.  11671017 and 11731009).\hspace{1mm} } \hspace{1mm}\hspace{1mm}
Renming Song\thanks{Research supported in part by a grant from the Simons
Foundation (\#429343, Renming Song).} \hspace{1mm}\hspace{1mm} and \hspace{1mm}\hspace{1mm}
Rui Zhang\footnote{Research supported in part by NSFC (Grant No. 11601354) and Youth Innovative Research Team of Capital Normal University. Corresponding author.}
\hspace{1mm} }
\date{}
\maketitle

\begin{abstract}
Suppose that $X=\{X_t, t\ge 0; \mathbb{P}_{\mu}\}$ is a supercritical superprocess in
a locally compact separable metric space $E$.  Let $\phi_0$ be a positive
 eigenfunction corresponding to the first eigenvalue $\lambda_0$ of
the generator of the mean semigroup of $X$. Then
$M_t:=e^{-\lambda_0t}\langle\phi_0, X_t\rangle$ is a positive martingale.
Let $M_\infty$ be the limit of $M_t$. It is known
(see, J. Appl. Probab. \textbf{46} (2009), 479--496) that
$M_\infty$ is non-degenerate iff  the $L\log L$ condition is  satisfied.
In this paper
we are mainly interested in the case when the $L\log L$ condition is not satisfied.
We prove that, under some conditions,  there exist
a positive function $\gamma_t$ on $[0, \infty)$
and  a non-degenerate random variable
$W$ such that for any
finite nonzero Borel measure $\mu$ on $E$,
\begin{equation*}
\lim_{t\to\infty}\gamma_t\langle \phi_0,X_t\rangle
=W,\qquad\mbox{a.s.-}\mathbb{P}_{\mu}.
\end{equation*}
We also give the almost sure limit of $\gamma_t\langle f,X_t\rangle$
for a class of general test functions $f$.
\end{abstract}

\medskip


\medskip


\maketitle

\bigskip

\baselineskip=6.0mm

\section{Introduction}

Suppose $\{Z_n, n\ge 0\}$ is a supercritical Galton-Watson process
with offspring number $L$. Let $m:=EL\in (1, \infty)$ be the mean of $L$.
Then $M_n:=\frac{Z_n}{m^n}$ is a non-negative martingale and thus has a finite limit $M_\infty$.
The well-known Kesten-Stigum theorem
says that the following three statements are equivalent:
  (i) $P(M_\infty=0)=P(Z_n=0 \mbox{ for large }n)$, i.\/e., the events $\{M_\infty=0\}$
and $\{ Z_n=0 \mbox{ for large }n\}$ are almost the same;
(ii) $EM_\infty=1$; (iii) $E(L\log L)<\infty$.
For a classical proof of this result, the reader is referred to the book \cite{AN72}
of Athreya and Ney.
In 1995, Lyons, Pemantle and  Peres \cite{LPP} gave
a  probabilistic proof of the
above $L\log L$ criterion of Kesten and Stigum.
If $E(L\log L)=\infty$, then $\lim_{n\to\infty}\frac{Z_n}{m^n}=0$ almost surely, which says that $m^{n}$ does not give the right growth
rate of  $Z_n$ conditional on non-extinction. It is natural to ask what is the right growth
rate of  $Z_n$.
In 1968, Seneta \cite{Sen68} proved that there is a sequence of positive numbers $c_n$
such that $c_nZ_n$ converges
to a non-degenerate random variable $W$ in distribution.
Heyde \cite{{Heyde70}} strengthened the convergence in distribution to
almost sure convergence.
 Later the problem of finding $c_n$ such that $c_nZ_n$
 converges to a non-degenerate limit is called the Seneta-Heyde norming problem.

Hoppe \cite{Hoppe} generalized the result of Heyde \cite{Heyde70} to
supercritical multitype branching processes,  Grey \cite{Grey74} proved a similar result for continuous state branching processes and Hering \cite{Hering} obtained  a similar result for supercritical
branching diffusions. In this paper we are going to consider the  Seneta-Heyde norming problem for general superprocesses under some conditions which are easy to check and satisfied by many superprocesses, including superdiffusions in a bounded domain and
also superprocesses with discontinuous spatial  motions.
We emphasize that we are mainly interested in the case when the $L\log L$ condition fails,
since the norming problem is already solved in  \cite{LRS09} for superdiffusions and  \cite{RSY16} for more general superprocesses when the $L\log L$ condition holds.

\subsection{Superprocesses and assumptions}

In this subsection,
we describe the setup of this paper and formulate our assumptions.

Suppose that $E$ is a locally compact separable metric space.
We will use
$\mathcal{B}(E)$ ($\mathcal{B}^+(E)$) to denote the
family of (non-negative)  Borel  functions on $E$,
$\mathcal{B}_b(E)$ ($\mathcal{B}_b^+(E)$) to denote the
family of (non-negative) bounded Borel  functions on $E$,
and $C(E)$($C_0(E)$, respectively)  to denote
the family of
continuous functions (vanishing at infinity, respectively) on $E$.

Suppose that $\partial$ is a separate point not contained
in $E$.
We will use $E_{\partial}$ to denote $E\cup\{\partial\}$.
Every function $f$ on $E$ is automatically extended to $E_{\partial}$ by setting $f(\partial)=0$.
 We will assume that $\xi=\{\xi_t,\Pi_x\}$ is a Hunt process on $E$ and $\zeta:=
 \inf\{t>0: \xi_t=\partial\}$ is the lifetime of $\xi$.
We will use $\{P_t:t\geq 0\}$ to denote the semigroup of $\xi$.
Suppose that $m$ is a
$\sigma$-finite Borel measure on $E$ with full support.
We will assume below that $\{P_t:t\geq 0\}$  has
a dual with respect to the
measure $m$ and the dual semigroup is
sub-Markovian.

The superprocess $X=\{X_t:t\ge 0\}$ we are going to work with is determined by two parameters:
a spatial motion $\xi=\{\xi_t, \Pi_x\}$ on $E$
which is a Hunt process,
and a branching mechanism $\varphi$ of the form
\begin{equation}\label{e:branm}
\varphi(x,s)=-\alpha(x)s+\beta(x)s^2+\int_{(0,+\infty)}(e^{-s \theta}-1+s \theta)n(x,{\rm d}\theta),
\quad x\in E, \, s\ge 0,
\end{equation}
where $\alpha\in \mathcal{B}_b(E)$, $\beta\in \mathcal{B}_b^+(E)$ and $n$ is a kernel from $E$ to $(0,\infty)$ satisfying
\begin{equation}\label{n:condition}
  \sup_{x\in E}\int_{(0,+\infty)} (\theta\wedge \theta^2) n(x,{\rm d}\theta)<\infty.
\end{equation}
Then there exists $M>0$ such that
\begin{equation*}
|\alpha(x)|+\beta(x)+
\int_{(0,+\infty)} (\theta\wedge \theta^2) n(x,{\rm d}\theta)\le M.
\end{equation*}
In this paper, we will exclude the case when
$\beta(\cdot)+n(\cdot,(0,\infty))=0$ $m$-almost everywhere.

Let ${\cal M}_F(E)$
be the space of finite measures on $E$, equipped with the topology of weak convergence.
The superprocess $X$  with spatial motion $\xi$ and branching mechanism $\varphi$
is a  Markov process taking values in ${\cal M}_F(E)$.
The existence of such superprocesses is well-known, see, for instance,
\cite{Dawson}, \cite{E.B.} or \cite{Li11}.
For any $\mu \in \mathcal{M}_F(E)$, we denote the
law
of $X$ with initial configuration $\mu$ by $\P_\mu$.
As usual,
$\langle f,\mu\rangle:=\int_E f(x)\mu(dx)$
and $\|\mu\|:=\langle 1,\mu\rangle$.
Throughout this paper, a real-valued function $u(t,x)$ on $[0,\infty)\times E_\partial$
is said to be locally bounded if, for any $t>0$, $\sup_{s\in[0,t], x\in E_\partial}|u(s, x)|<\infty$.
According to \cite[Theorem 5.12]{Li11}, there is a Hunt process
$X=\{\Omega, {\cal G}, {\cal G}_t, X_t, \P_\mu\}$ taking values in  $\mathcal{M}_F(E)$
such that for every
$f\in \mathcal{B}^+_b(E)$ and $\mu \in \mathcal{M}_F(E)$,
\begin{equation}
  -\log \P_\mu\left(e^{-\langle f,X_t\rangle}\right)
  =\langle V_tf,\mu\rangle,
\end{equation}
where $V_tf(x)$
is the unique locally bounded non-negative solution to the equation
\begin{equation}\label{u_f}
   V_tf(x)+\Pi_x\int_0^t\varphi(\xi_s, V_{t-s}f(\xi_s)){\rm d}s=\Pi_x f(\xi_t), \quad x\in E_\partial,
\end{equation}
where we use the convention that
$\varphi(\partial,s)=0$ for all $s\ge0$.
Since $f(\partial)=0$, we have
$V_tf(\partial)=0$ for any $t\ge 0$.
In this paper, the superprocess we deal with is always this Hunt realization.

For any $f\in\mathcal{B}_b(E)$ and $(t, x)\in (0, \infty)\times E$, we define
\begin{equation}\label{1.26}
   T_tf(x):=\Pi_x \left[e^{\int_0^t\alpha(\xi_s){\rm d}s}f(\xi_t)\right].
\end{equation}
It is well known that $T_tf(x)=\P_{\delta_x}\langle f,X_t\rangle$ for every $x\in E$.

We will always assume that there exists a family of
continuous and strictly positive functions $\{p(t,x,y):t>0\}$ on $E\times E$ such that for any $t>0$ and nonnegative function $f$ on $E$,
\begin{equation*}
  P_tf(x)=\int_E p(t,x,y)f(y)m({\rm d}y).
\end{equation*}
Define
\begin{equation*}
a_t(x):=\int_E p(t,x,y)^2\,m({\rm d}y),\qquad \hat{a}_t(x):=\int_E p(t,y,x)^2\,m({\rm d}y).
\end{equation*}
Our first main assumption is
\begin{assumption}\label{assum1}
\begin{description}
  \item[(i)]  For any $t>0$, $\int_E p(t,x,y)\,m({\rm d}x)\le 1$.
  \item[(ii)] For any $t>0$, we have
  \begin{equation}\label{1.17}
     \int_E a_t(x)\,m({\rm d}x)=\int_E\hat{a}_t(x)\,
m({\rm d}x)=\int_E\int_E p(t,x,y)^2\,m({\rm d}y)\,m({\rm d}x)<\infty.
  \end{equation}
  Moreover, the functions $x\to a_t(x)$ and $x\to\hat{a}_t(x)$ are continuous on $E$.
\end{description}
\end{assumption}

Note that, in Assumption \ref{assum1}(i), the integration is with respect to the first
space variable. It implies that the dual semigroup $\{{\widehat P}_t:t\ge 0\}$
of $\{P_t:t\ge 0\}$ with respect to $m$ defined by
\begin{equation*}
{\widehat P}_tf(x)=\int_Ep(t,y, x)f(y)m({\rm d}y)
\end{equation*}
is sub-Markovian. Assumption \ref{assum1}(ii) is a pretty weak $L^2$ condition and it allows
us to apply results on operator semigroups in Hilbert spaces.
By H\"older's inequality, we have
\begin{equation}\label{1.1}
 p(t+s,x,y)=\int_E p(t,x,z)p(s,z,y)\,m({\rm d}z)\le (a_t(x))^{1/2}(\hat{a}_s(y))^{1/2}.
\end{equation}

It is well known and easy to check that $\{P_t:t\ge 0\}$
and $\{\widehat{P}_t:t\ge 0\}$ are strongly continuous
contraction semigroups on $L^2(E, m)$,
see \cite{RSZ3} for a proof.
We will use $\langle\cdot, \cdot\rangle_m$ to denote the inner product in $L^2(E, m)$.
Since $p(t, x, y)$ is continuous in $(x, y)$,
by \eqref{1.1}, Assumption \ref{assum1}(ii) and
the dominated convergence theorem, we have that,
for any $f\in L^2(E,m)$,  $P_tf$  and $\widehat{P}_tf$ are continuous.

It follows from Assumption \ref{assum1}(ii) that,
for each $t>0$, $P_t$ and $\widehat{P}_t$
are compact operators on $L^2(E,m)$.
Let $\widetilde{L}$ and $\widehat{\widetilde{L}}$ be the infinitesimal generators of the semigroups $\{P_t\}$ and $\{\widehat{P}_t\}$ in $L^2(E,m)$ respectively.
Define $\widetilde{\lambda}_0:=\sup \Re(\sigma(\widetilde{L}))=\sup\Re(\sigma(\widehat{\widetilde{L}}))$.
By Jentzsch's theorem (Theorem V.6.6 on page 337 of \cite{Sch}),
$\widetilde{\lambda}_0$ is an eigenvalue of multiplicity 1 for both $\widetilde{L}$ and
$\widehat{\widetilde{L}}$, and that an eigenfunction $\widetilde{\phi}_0$ of $\widetilde{L}$ corresponding
to $\widetilde{\lambda}_0$ can be chosen to be strictly positive $m$-almost
everywhere with $\|\widetilde{\phi}_0\|_2=1$ and
an eigenfunction $\widetilde{\psi}_0$ of $\widehat{\widetilde{L}}$ corresponding to
$\widetilde{\lambda}_0$  can be chosen to be strictly positive $m$-almost
everywhere with $\langle \widetilde{\phi}_0, \widetilde{\psi}_0\rangle_m=1$.
Thus for $m$-almost every $x\in E$,
\begin{equation*}
e^{\widetilde{\lambda}_0}\widetilde{\phi}_0(x)=P_1\widetilde{\phi}_0(x),
\qquad
e^{\widetilde{\lambda}_0}\widetilde{\psi}_0(x)=\widehat{P}_1\widetilde{\psi}_0(x).
\end{equation*}
Hence $\widetilde{\phi}_0$ and $\widetilde{\psi}_0$ can be chosen to be continuous and strictly
positive everywhere on $E$.

Our second assumption is

\begin{assumption}\label{assum3}
\begin{description}
 \item[(i)] $\widetilde{\phi}_0$ is bounded.

 \item[(ii)]
The semigroup $\{P_t,t\ge0\}$ is intrinsically ultracontractive,
that is, there exists $c_t>0$ such that
\begin{equation}\label{newcondition2}
  p(t,x,y)\le c_t\widetilde{\phi}_0(x)\widetilde{\psi}_0(y).
\end{equation}
\end{description}
\end{assumption}

Assumption \ref{assum3} is a pretty strong assumption on the semigroup
$\{P_t:t\ge0\}$. However, this assumption is satisfied in a lot of cases.
See \cite[Section 1.4]{RSZ4} for examples of  Markov processes
satisfying the assumption above.
The concept of intrinsic ultracontractivity was introduced by Davies and Simon
in \cite{DS} in the setting of symmetric semigroups. This concept was extended to
the non-symmetric setting in \cite{KiSo08b, KiSo08c, KiSo09}.
Intrinsic ultracontractivity has been studied intensively in the last 30 years
and there are many results on the intrinsic ultracontractivity of
semigroups, see \cite{KiSo08b, KiSo08c, KiSo09} and the references therein.

We have proved in \cite[Lemma 2.1]{RSZ4} that
there exists a function $q(t, x, y)$ on $(0, \infty)\times E\times E$ which is continuous
in $(x, y)$ for each $t>0$ such that
\begin{equation}\label{comp0}
e^{-Mt}p(t,x,y) \le q(t,x,y)\le e^{Mt}p(t,x,y), \quad (t, x, y)\in (0, \infty)\times E\times E,
\end{equation}
and that for any bounded Borel function $f$ and any $(t, x)\in (0, \infty)\times E$,
\begin{equation*}
  T_tf(x)=\int_E q(t,x,y)f(y)\,m({\rm d}y).
\end{equation*}
It follows immediately that
\begin{equation}\label{Lp}
  \|T_tf\|_2\le e^{Mt}\|P_tf\|_2\le e^{Mt}\|f\|_2.
\end{equation}
In \cite{RSZ3}, we have proved that $\{T_t: t\ge 0\}$ is a strongly continuous semigroup on $L^2(E, m)$.
Let $\{\widehat{T}_t,t>0\}$ be the adjoint semigroup
on $L^2(E,m)$ of $\{T_t,t>0\}$, that is,
for $f\in L^2(E,m)$,
\begin{equation*}
\widehat{T}_tf(x)=\int_E q(t,y,x)f(y)\,m({\rm d}y).
\end{equation*}
We have proved in \cite{RSZ3} that $\{\widehat{T}_t: t\ge 0\}$ is also a strongly continuous semigroup on $L^2(E, m)$.
We claim that, for all $t>0$ and $f\in L^2(E,m)$, $T_tf$  and $\widehat{T}_tf$ are continuous.
In fact, since $q(t, x, y)$ is continuous in $(x, y)$, by \eqref{1.1}, \eqref{comp0},
Assumption 1(ii)
and the dominated convergence theorem,
we have that, for any $f\in L^2(E,m)$,  $T_tf$  and $\widehat{T}_tf$ are continuous.

By Assumption \ref{assum1}(ii) and \eqref{comp0}, we get that
\begin{equation*}
\int_E\int_E q^2(t,x,y)\,
m({\rm d}x)\,m({\rm d}y)\le e^{2Mt}\int_E\int_E p^2(t,x,y)\,m({\rm d}x)\,m({\rm d}y)
<\infty.
\end{equation*}
Thus, for each $t>0$, $T_t$ and $\widehat{T}_t$
are compact operators on $L^2(E,m)$.
Let $L$ and $\widehat{L}$ be the infinitesimal generators of the semigroups $\{T_t\}$ and $\{\widehat{T}_t\}$ in $L^2(E,m)$ respectively.
Define $\lambda_0:=\sup \Re(\sigma(L))=\sup\Re(\sigma(\widehat{L}))$.
By Jentzsch's theorem,
$\lambda_0$ is an eigenvalue of multiplicity 1 for both $L$ and
$\widehat{L}$, and that an eigenfunction $\phi_0$ of $L$ corresponding
to $\lambda_0$ can be chosen to be strictly positive $m$-almost
everywhere with $\|\phi_0\|_2=1$ and
an eigenfunction $\psi_0$ of $\widehat{L}$ corresponding to
$\lambda_0$  can be chosen to be strictly positive $m$-almost
everywhere with $\langle \phi_0,\psi_0\rangle_m=1$.
Thus for $m$-almost every $x\in E$,
\begin{equation*}
e^{\lambda_0}\phi_0(x)=T_1\phi_0(x),
\qquad
e^{\lambda_0}\psi_0(x)=\widehat{T}_1\psi_0(x).
\end{equation*}
Hence $\psi_0$ and $\phi_0$ can be chosen to be continuous and strictly
positive everywhere on $E$.

Using Assumption \ref{assum3}, the boundedness of $\alpha$ and an
argument similar to that used in the proof of \cite[Theorem 3.4]{DS},
one can show the following:
\begin{description}
\item[(i)] $\phi_0$ is bounded.

\item[(ii)]
The semigroup $\{T_t,t\ge0\}$ is intrinsically ultracontractive,
that is, there exists $c_t>0$ such that
\begin{equation}\label{condition2}
  q(t,x,y)\le c_t\phi_0(x)\psi_0(y).
\end{equation}
\end{description}

The main interest of this paper is on supercritical superprocesses, so we assume that
\begin{assumption}\label{assum2}
  $\lambda_0>0$.
\end{assumption}

Define $q_t(x):=\P_{\delta_x}(\|X_t\|=0)$.
Since $\P_{\delta_x}\|X_t\|=T_t1(x)>0$, we have $q_t(x)<1$.
Note that $q_t(x)$ is non-decreasing in $t$. Hence  the limit
\begin{equation*}
q(x):=\lim_{t\to\infty}q_t(x)
\end{equation*}
exists.
 Then $q(x)=\P_{\delta_x}\{\|X_t\|=0 \mbox{ for some } t>0\}$ is the extinction probability.
In this paper, we also assume that
\begin{assumption}\label{assum4}
There exists $t_0>0$ such that
\begin{equation}\label{4.2'}
  \inf_{x\in E} q_{t_0}(x)>0.
\end{equation}
\end{assumption}

In \cite[Section 2.2]{RSZ4}, we gave
a sufficient condition (in term of the branching mechnism $\varphi$) for Assumption \ref{assum4}. In particularly, if $\inf_{x\in E}\beta(x)>0$, then Assumption \ref{assum4} holds.
In Lemma \ref{lem:ext}, we will show that, under our assumptions,
$q(x)<1,$ for all $x\in E.$

\subsection{Main results}

Define
$$M_t:=e^{-\lambda_{0}t}\langle \phi_0,X_{t}\rangle,\quad t\ge 0.$$
It is easy to prove that (see \cite[Theorem 3.2]{RSY16}, for example),
for every $\mu\in\mathcal{M}_F(E)$,
$\{M_t, t\ge 0\}$
is a non-negative $\mathbb{P}_{\mu}$-martingale with respect to the filtration $\{\mathcal{G}_{t}:t\ge 0\}$.
Thus $\{M_t, t\ge 0\}$ has a $\mathbb{P}_{\mu}$-a.s. finite limit denoted as $M_\infty$.

Let
$n^{\phi_0}(x,{\rm d}\theta)$
be the kernel from $E$ to $(0,\infty)$ defined by
\begin{equation*}
\int_0^\infty f(\theta)n^{\phi_0}(x,{\rm d}\theta)=\int_0^\infty f(\theta\phi_0(x))n(x,{\rm d}\theta).
\end{equation*}
By the boundedness of $\phi_0$ and the assumption \eqref{n:condition}, we get that
there exists $\widetilde{M}>0$ such that
\begin{equation}\label{n:condition3}
\sup_{x\in E}\int_0^\infty (\theta\wedge \theta^2)n^{\phi_0}(x,{\rm d}\theta)\le \widetilde{M}.
\end{equation}
Let
$l(x):=\int_1^\infty \theta\log \theta \,n^{\phi_0}(x, {\rm d}\theta).$
The following $L\log L$ criterion was proved for superdiffusions in \cite{LRS09}.

\emph{
\textbf{$L\log L$ criterion:}
$M_\infty$ is non-degenerate under $\mathbb P_\mu$ for all
 nonzero finite measures $\mu$ on $E$ if and only if
\begin{equation}\label{llog2}
\int_E\psi_0(x)l(x)m({\rm d}x)<\infty.
\end{equation}}

At first glance, the roles of $\phi_0$ and $\psi_0$ are not symmetric in \eqref{llog2}. This is not the case. In fact,
\eqref{llog2} is equivalent to
\begin{equation}\label{llog2'}
\int_E\phi_0(x)\psi_0(x)m({\rm d}x)\int^\infty_1(\theta\log \theta)n(x, {\rm d}\theta)<\infty,
\end{equation}
which says that the spatial average of the
``$\theta\log \theta$"
moment $\int^\infty_1(\theta\log \theta)n(x, {\rm d}\theta)$
with respect to  the probability measure $\phi_0(x)\psi_0(x)m({\rm d}x)$ is finite.
Note that
\begin{align*}
  l(x)&=\int_1^\infty \theta\log \theta\,n^{\phi_0}(x,{\rm d}\theta)= \phi_0(x)\int_{\phi_0(x)^{-1}}^\infty \theta(\log \theta+\log\phi_0(x))\,n(x,{\rm d}\theta)\\
  &=\phi_0(x)\int_{1}^\infty \theta\log \theta\,n(x,{\rm d}\theta)+\phi_0(x)\int_{\phi_0(x)^{-1}}^1 \theta\log \theta\,n(x,{\rm d}\theta)\\
  &\ \ +\phi_0(x)\log\phi_0(x)\int_{\phi_0(x)^{-1}}^\infty \theta\,n(x,{\rm d}\theta).
\end{align*}
Since $\phi_0$ is bounded,  $\phi_0(x)|\log\phi_0(x)|$ is bounded above, say by $C$. Thus,
$$\sup_{x\in E}\phi_0(x)|\log\phi_0(x)|\int_{\phi_0(x)^{-1}}^\infty \theta\,n(x,{\rm d}\theta)\le C\sup_{x\in E}\int_{\|\phi_0\|_\infty^{-1}}^\infty \theta\,n(x,{\rm d}\theta)<\infty.$$
Note also that
\begin{align*}
  &\sup_{x\in E}\phi_0(x)\left|\int_{\phi_0(x)^{-1}}^1 \theta\log \theta\,n(x,{\rm d}\theta)\right|\le\sup_{x\in E}\phi_0(x)|\log\phi_0(x)|\left|\int_{\phi_0(x)^{-1}}^1\theta\,n(x,{\rm d}\theta)\right|\\
  &\le C\sup_{x\in E}\Big(\textbf{1}_{\phi_0(x)>1}\int_{\|\phi_0\|_\infty^{-1}}^1 \theta\,n(x,{\rm d}\theta)+\textbf{1}_{\phi_0(x)\le1}\int_1^\infty \theta\,n(x,{\rm d}\theta)\Big)<\infty.
\end{align*}
In Section 2, we will show that $\psi_0\in L^1(E,m)$ (see the first paragraph of Section 2). Thus \eqref{llog2} is equivalent to \eqref{llog2'}.

Recently, the $L\log L$ criterion above  was extended to more general
superprocesses with possible non-local branching mechanisms in \cite{RSY16}.

The $L\log L$ criterion above says that, under condition \eqref{llog2},
$e^{\lambda_0t}$
gives the growth rate of $\langle\phi_0, X_t\rangle$ as $t\to\infty$ conditioned
on non-extinction.
However, when condition \eqref{llog2} is not satisfied, the theorem above
does not provide much information about the growth rate of
$\langle\phi_0, X_t\rangle$.

The first objective of this paper is to solve the Seneta-Heyde norming
problem for the martingale $M_t$, that is,
to find a positive function $\gamma_t$  on $[0, \infty)$
such that $\gamma_t\langle\phi_0, X_t\rangle$ has a non-degenerate limit
as $t\to\infty$.
Although our results (Theorems \ref{main1}, \ref{LLN} and \ref{LLN2}) also cover the case when \eqref{llog2} holds, only the results
in the case when \eqref{llog2} fails are new, see Theorem \ref{llogl3} below.
It is easy to find examples such that \eqref{llog2} fails. For example, if
$n(x, {\rm d}\theta)=c(x)[I_{(0, 2)}(\theta)+I_{[2, \infty)}(\theta)\theta^{-2}(\log \theta)^{-\beta}] {\rm d}\theta$
with $\beta\in(1, 2]$ and $c(x)$ being  strictly positive bounded measurable function on $E$, then \eqref{n:condition} is satisfied, but \eqref{llog2'}, which is equivalent to \eqref{llog2}, fails.

Let $v(x):=-\log q(x)$. By the branching property of $X$,  we have
\begin{equation*}
 \P_{\mu}(\|X_t\|=0 \mbox{ for some } t>0)=e^{-\la v,\mu\ra}.
\end{equation*}
\begin{theorem}\label{main1}
There exist a positive function $\gamma_t$ on $[0, \infty)$ and
a non-degenerate random variable $W$ such that
$$\lim_{t\to\infty}\frac{\gamma_{t}}{\gamma_{t+s}}=e^{\lambda_0s},\quad\forall s\ge 0,$$
and that for any nonzero $\mu\in\mathcal{M}_F(E)$,
\begin{equation*}
\lim_{t\to\infty}\gamma_t\langle \phi_0,X_t\rangle
=W,\qquad
\mbox{a.s.-}\P_{\mu}
\end{equation*}
and
\begin{equation*}
\P_{\mu}(W=0)=e^{-\la v,\mu\ra},\qquad \P_{\mu}(W<\infty)=1.
\end{equation*}
\end{theorem}

Moreover, we have the following $L\log L$ criterion.

\begin{theorem}\label{llogl3}
The following conditions are equivalent:
\begin{description}
\item[(1)] $M_\infty$ is non-degenerate for some nonzero $\mu\in{\cal M}_F(E)$;
  \item[(2)] $M_\infty$ is non-degenerate for all nonzero $\mu\in{\cal M}_F(E)$;
  \item[(3)] $l_0:=\lim_{t\to\infty}e^{\lambda_0t}\gamma_t\in (0, \infty)$;
  \item[(4)] $\int_E\psi_0(x)l(x)m({\rm d}x)<\infty;$
 \item[(5)] $\P_{\mu}W<\infty$ for some nonzero $\mu\in{\cal M}_F(E)$.
 \item[(6)] $\P_{\mu}W<\infty$ for all  $\mu\in{\cal M}_F(E)$.
\end{description}
\end{theorem}

Further properties of the limit random variable $W$, such as absolute continuity and tail probabilities, are studied in \cite{RSZ19}, a sequel to the present paper.

The second objective of this paper is to study
the almost sure limit behaviour of $\gamma_t\langle f, X_t\rangle$ as $t\to\infty$ for
a class of bounded continuous functions $f$.
It turns out that, for $f$ belonging to this class,
$\lim_{t\to\infty}\gamma_t\langle f, X_t\rangle=\la f,\psi_0\ra_m W$,
$\mathbb{P}_\mu$-a.s.
for any nonzero $\mu\in{\cal M}_F(E)$,
see Theorems \ref{LLN} and \ref{LLN2}.

The rest of the paper is organized as follows.
Section 2 contains our basic estimates and Section 3 deals with some properties of the
extinction probability.
In Section 4, we will  define and investigate backward iterates, which is needed in the proof of Theorem \ref{main1}.  The proofs of Theorem \ref{main1} and  Theorem \ref{llogl3} are  given in Section 5.  We remark that  we will prove  Theorem \ref{llogl3}  without using
the $L\log L$ criterion in \cite{LRS09}.
The strong limit behaviour of $\gamma_t\langle f, X_t\rangle$ as $t\to\infty$ for a class of general bounded continuous functions $f$ is given in Section 6.
In the last section, we give some concluding remarks.

In the remainder of this paper, $C$ will stand for a constant whose value might change from one appearance to the next.

\section{Some estimates}

According to \cite[Thorem 2.7]{KiSo08b}, under Assumptions \ref{assum1}--\ref{assum3},
for any $\delta>0$, there exist constants $\gamma=\gamma(\delta)>0$ and $c=c(\delta)>0$ such that,
for any $(t,x,y)\in [\delta,\infty)\times E\times E$, we have
\begin{equation}\label{density0}
  \left|e^{-\lambda_0t}q(t,x,y)-\phi_0(x)\psi_0(y)\right|\le ce^{-\gamma t}\phi_0(x)\psi_0(y).
\end{equation}
Take $t$ large enough so that $ce^{-\gamma t}<\frac12$, then we have
\begin{equation*}
e^{-\lambda_0 t}q(t, x, y)\ge \frac12\phi_0(x)\psi_0(y).
\end{equation*}
Since $q(t, x, \cdot)\in L^1(E, m)$, we have $\psi_0\in L^1(E, m)$.

It follows from \eqref{density0} that, if $f\in\mathcal{B}_b^+(E)$
 then  $\la f,\psi_0\ra_m<\infty$ and  for any $(t,x)\in[\delta,\infty)\times E$,
\begin{equation}\label{T_t}
  \left|e^{-\lambda_0t}T_tf(x)-\langle f,\psi_0\rangle_m \phi_0(x)\right|
  \le ce^{-\gamma t}\langle |f|,\psi_0\rangle_m \phi_0(x)
\end{equation}
 and
 \begin{equation}\label{1.6}
   (1-ce^{-\gamma t})\langle |f|,\psi_0\rangle_m\phi_0(x)\le e^{-\lambda_0t}T_t|f|(x)\le  (1+c)\langle |f|,\psi_0\rangle_m\phi_0(x).
 \end{equation}
For $x\in E$ and $s>0$, we define
\begin{equation*}\label{def-r}
r(x,s)=\varphi(x,s)+\alpha(x)s.
\end{equation*}

\begin{lemma}\label{lem:1.3}
For any $H\ge 1$,
\begin{equation}\label{2.6}
0\le r(x,s)\le (3/2+H/2)Ms^2+\int_H^\infty \theta\,n(x,{\rm d}\theta)s.
\end{equation}
\end{lemma}
\noindent\textbf{Proof:}
By definition,
\begin{equation}\label{2.2}
r(x,s)=\beta(x)s^2+\int_0^\infty(e^{-\theta s}-1+\theta s)\,n(x,{\rm d}\theta).
\end{equation}
Note that for any $\theta>0$,
\begin{equation}\label{2.7}
  0<e^{-\theta}-1+\theta\le
  \left(\frac\theta2\wedge1\right)\theta.
\end{equation}
Thus for any $H\ge1$,
\begin{align*}
  r(x,s)&\le Ms^2+\frac{1}{2}s^2\int_0^H\theta^2\,n(x,{\rm d}\theta)+s\int_H^\infty \theta\,n(x,{\rm d}\theta) \\
   &\le Ms^2+\frac{1}{2}s^2\left[\int_0^1\theta^2\,n(x,{\rm d}\theta)
   +H\int_1^H \theta\,n(x,{\rm d}\theta)\right]+s\int_H^\infty \theta\,n(x,{\rm d}\theta) \\
   &\le (3/2+H/2)Ms^2+\int_H^\infty \theta\,n(x,{\rm d}\theta)s.
\end{align*}
\hfill $\Box$

For any $f\in \mathcal{B}_b^+(E)$ satisfying $m(f>0)>0$,
define, for any $(t, x)\in [0, \infty)\times E$,
\begin{equation}\label{de:R1}
 R_f(t,x):=T_tf(x)-V_tf(x),
\end{equation}
and
\begin{equation}\label{1.10}
  g_f(t,x):=\frac{R_f(t,x)}{e^{\lambda_0 t}\langle f,\psi_0\rangle_m\phi_0(x)}.
\end{equation}

\begin{lemma}\label{lem1.2}
Assume that $f\in \mathcal{B}_b^+(E)$ and $m(f>0)>0$. Then
\begin{align}\label{1.9}
 R_f(t,x)\ge 0,\quad (t, x)\in [0, \infty)\times E,
\end{align}
and
\begin{equation}\label{1.22}
  \lim_{\|f\|_\infty\to0}\|g_f(t)\|_\infty=0,\quad t\ge0.
\end{equation}
\end{lemma}

\noindent\textbf{Proof:}
It follows from \cite[Theorem 2.23]{Li11} and \eqref{u_f} that $V_tf(x)$ also satisfies
\begin{equation}\label{1.8}
  V_tf(x)=-\int_0^t T_s\left[r(\cdot,V_{t-s}f(\cdot))\right](x){\rm d}s+T_tf(x), \quad t\ge 0, x\in E,
\end{equation}
which implies that
\begin{equation}\label{IF}
R_f(t,x)=\int_0^t T_s\left[r(\cdot,V_{t-s}f(\cdot)\right](x){\rm d}s,\quad t\ge 0, x\in E.
\end{equation}
Since $r(x,s)\ge 0$, we have $R_f(t,x)\ge 0$,
i.\/e., (\ref{1.9}) holds.

Since $T_0f(x)=V_0f(x)=f(x)$, then $R_f(0,x)=0$, which implies that $g_f(0,x)=0$.
Thus, it suffices to prove that, for any $\delta>0$,  \eqref{1.22} is true for $t>\delta$.
It follows from Lemma \ref{lem:1.3} that for any $H\ge 1$,
\begin{align*}
  r(x,V_{t-s}f(x))&\le (2+H)M(V_{t-s}f(x))^2+V_{t-s}f(x)\int_H^\infty \theta\,n(x,{\rm d}\theta) \\
 &\le (2+H)M(T_{t-s}f(x))^2+T_{t-s}f(x)
 \int_H^\infty \theta\,n(x,{\rm d}\theta).
\end{align*}
By \eqref{IF}, we have
\begin{align*}
  R_f(t,x)&\le (2+H)M\int_0^t T_s[(T_{t-s}f)^2](x){\rm d}s+\int_0^tT_s\left[T_{t-s}f
  \int_H^\infty \theta\,n(\cdot,{\rm d}\theta)\right](x){\rm d}s\\
 &=:(I)+(II).
  \end{align*}

For part (I), since $T_{t-s}f(x)\le \|f\|_\infty e^{Mt}$, we have $T_s[(T_{t-s}f)^2](x)\le e^{Mt}\|f\|_\infty T_tf(x)$. Thus, by \eqref{1.6}, we have that, for any $t>\delta$,
\begin{align*}
  (I) &\le  (2+H)Me^{Mt}t\|f\|_\infty T_tf(x)\nonumber\\
   &\le  (2+H)Me^{Mt}t (1+c(\delta))e^{\lambda_0t}\langle f, \psi_0\rangle_m\phi_0(x)\|f\|_\infty.
\end{align*}
Hence we have that, for any $t>\delta$,
\begin{equation}\label{EI}
  \lim_{\|f\|_\infty\to0}\sup_{x\in E}\frac{(I)}{e^{\lambda_0t}\langle f, \psi_0\rangle_m\phi_0(x)}=0.
\end{equation}

Now we deal with part (II). For any $H>1$, $t>\delta$ and  $0< \epsilon<t$,
\begin{equation}\label{EII1}
 \left(\int_0^{\epsilon}+\int_{t-\epsilon}^t\right) T_s\Big[T_{t-s}f\int_H^\infty
\theta \,n(\cdot,{\rm d}\theta)\Big](x){\rm d}s\le
 2\epsilon MT_tf(x)
 \le 2\epsilon M(1+c(\delta))e^{\lambda_0t}\langle f, \psi_0\rangle_m\phi_0(x).
\end{equation}
For any $\epsilon<s<t-\epsilon$, by \eqref{1.6},
$$
T_{t-s}f(x)\le (1+c(\epsilon))e^{\lambda_0t}\langle f, \psi_0\rangle_m\phi_0(x).
$$
Hence, we have
\begin{align}
  &\int_\epsilon^{t-\epsilon} T_s\Big[T_{t-s}f
\int_H^\infty \theta\,n(\cdot,{\rm d}\theta)\Big](x){\rm d}s\nonumber\\
  &\le (1+c(\epsilon))e^{\lambda_0t}\langle f, \psi_0\rangle_m\int_\epsilon^{t-\epsilon}
  T_s\Big[\int_H^\infty \theta\,n(\cdot,{\rm d}\theta)\phi_0\Big](x){\rm d}s\nonumber\\
  &\le(1+c(\epsilon))^2e^{2\lambda_0t} t\Big\la
  \int_H^\infty \theta\,n(\cdot,{\rm d}\theta)\phi_0,\psi_0
  \Big\ra_m\langle f, \psi_0\rangle_m\phi_0(x).\label{EII2}
\end{align}
Thus, combining \eqref{EII1} and \eqref{EII2}, we get that
\begin{equation*}
  \limsup_{\|f\|_\infty\to0}\sup_{x\in E}\frac{(II)}{e^{\lambda_0t}
  \langle f, \psi_0\rangle_m\phi_0(x)}\le 2\epsilon M(1+c(\delta))+(1+c(\epsilon))^2e^{\lambda_0t} t\Big\la
 \int_H^\infty \theta\,n(\cdot,{\rm d}\theta)\phi_0,\psi_0\Big\ra_m.
\end{equation*}
Now, first letting $H\to\infty$ and then $\epsilon\to0$,
applying the monotone convergence theorem,
we get that
\begin{equation}\label{EII}
  \lim_{H\to\infty}\limsup_{\|f\|_\infty\to0}\sup_{x\in E}\frac{(II)}{e^{\lambda_0t}\langle f, \psi_0\rangle_m\phi_0(x)}=0.
\end{equation}
Combining \eqref{EI} and \eqref{EII}, we get that
\begin{equation*}
\lim_{\|f\|_\infty\to0}\sup_{x\in E}\frac{R_{f}(t,x)}{e^{\lambda_0t}\langle f, \psi_0\rangle_m\phi_0(x)}=0.
\end{equation*}
\hfill $\Box$

\section{Extinction probability}
Recall that, for any $t>0$ and $x\in E$,
$$
q_t(x)=\P_{\delta_x}(\|X_t\|=0)\qquad \mbox{and }\qquad q(x)=\lim_{t\to\infty}q_t(x).
$$

\begin{lemma}\label{lem:ext}
For any $x\in E$,
$$q(x)<1.$$
\end{lemma}
\noindent\textbf{Proof:}
Let $v_t(x):=-\log q_t(x)$. Recall that $v(x)=-\log q(x)$.
Since $q_t(x)<1$, we have $v_t(x)>0$.
By Assumption \ref{assum4}, we have
for $s>t_0$,
$$
\|v\|_\infty\le \|v_{s}\|_\infty\le \|v_{t_0}\|_\infty=-\log\left(\inf_{x\in E}q_{t_0}(x)\right)<\infty.
$$
Recall that, for $\theta>0$,
$$
V_t\theta(x)=-\log \P_{\delta_x}e^{-\langle \theta,X_t\rangle}.
$$
By the Markov property of $X$,
\begin{equation}\label{1.20'}
  q_{t+s}(x)=\lim_{\theta\to\infty}\P_{\delta_x}\left(e^{-\theta\|X_{t+s}\|}\right)
  =\lim_{\theta\to\infty}\P_{\delta_x}
  \left(e^{-\langle V_s\theta,X_t\rangle}\right)
  =\P_{\delta_x}\left(e^{-\langle v_s ,X_t\rangle}\right).
\end{equation}
It follows from \eqref{1.10} that, for any $s>t_0$,
\begin{equation}\label{2.4}
  v_{t+s}(x)=V_t(v_s)(x) = T_{t}(v_s)(x)-R_{v_s}(t,x)=T_{t}(v_s)(x)-g_{v_s}(t,x)e^{\lambda_0t}\langle v_s,\psi_0\rangle_m\phi_0(x).
\end{equation}
Thus, for $s>t_0$, we have
\begin{equation*}\label{v:eq}
  \langle v_{t+s},\psi_0\rangle_m\ge (1-\|g_{v_s}(t)\|_\infty)e^{\lambda_0t}\langle v_s,\psi_0\rangle_m.
\end{equation*}
Since $v_t(x)$ is positive and non-increasing in $t$,
we have that for all $t>0$ and $s>t_0$,
\begin{equation}\label{2.3}
  (1-\|g_{v_s}(t)\|_\infty)e^{\lambda_0t}\le 1.
\end{equation}

We claim that $\langle v,\psi_0\rangle_m>0$. Otherwise, $\langle v,\psi_0\rangle_m=0$.
By \eqref{2.4}, we have
$$
\|v_{1+s}\|_\infty\le \|T_{1}(v_s)\|_\infty\le (1+c)e^{\lambda_0}\langle v_s,\psi_0\rangle_m\|\phi_0\|_\infty\to0
$$
as $s\to\infty$. Thus $\lim_{s\to\infty}\|v_{s}\|_\infty=0$.
Hence by \eqref{1.22},
$$
\lim_{s\to\infty}\|g_{v_s}(t)\|_\infty=0.
$$
It follows from \eqref{2.3} that, for all $t>0$, $e^{\lambda_0t}\le 1$, which is a contradiction to the assumption that $\lambda_0>0$.
Hence the claim above is valid.

By letting $s\to\infty$ in \eqref{1.20'},  we get that
\begin{equation}\label{eq:v}
\exp\{-v(x)\}=q(x)=\P_{\delta_x}\exp\left\{-\la v,X_t\ra\right\}.
\end{equation}
Let $c$ and $\gamma$ be the constants in \eqref{T_t} with $\delta=1$.
For $t$ large enough, we
have $1-ce^{-\gamma t}>0$. Thus for $t$ large enough, we have
$$
 T_tv(x)\ge (1-ce^{-\gamma t})e^{\lambda_0t}\langle v,\psi_0\rangle_m\phi_0(x)>0.
$$
Hence for all $x\in E$ and $t$ large enough,
$$
\P_{\delta_x}(\langle v,X_t\rangle>0)>0,
$$
which implies that
$$
q(x)=\P_{\delta_x}(e^{-\langle v,X_t\rangle})<1.
$$
The proof is now complete.
\hfill $\Box$

\begin{lemma}\label{lem5}
$V:=\lim_{t\to\infty}\la v,X_t\ra\in[0,\infty]$ exists, and
satisfies that, for all $x\in E$,
$$
\P_{\delta_x}(V=0)=\exp\{-v(x)\}
=q(x)
$$
and
$$
\P_{\delta_x}(V =\infty)=1-\exp\{-v(x)\}
  =1-q(x).
$$
Moreover, for any $\theta>0$, we have
\begin{equation}\label{3.11}
  \lim_{t\to\infty}V_t(\theta v)(x)=v(x),\quad x\in E.
\end{equation}
\end{lemma}

\noindent\textbf{Proof:}
By \eqref{eq:v} and the Markov property of $X$,
$\{e^{-\la v,X_t\ra}, t\ge0\}$ is a bounded martingale.
Thus $\lim_{t\to\infty}\exp\left\{-\la v,X_t\ra\right\}$ exists and is in $[0, 1]$, which implies that $V:=\lim_{t\to\infty}\la v,X_t\ra\in[0,\infty]$ exists.

Since
$\exp\left\{-\la v,X_t\ra\right\}\le 1$, we have
\begin{equation}\label{3.1}
  \P_{\delta_x}\exp\{-V\}=\exp\{-v(x)\}.
\end{equation}
On the other hand,
\begin{align*}
\P_{\delta_x}(V=0)\ge \lim_{t\to\infty}\P_{\delta_x}(\|X_t\|=0)=\exp\{-v(x)\},
\end{align*}
which implies that
\begin{equation}\label{3.2}
  \P_{\delta_x}\exp\{-V\}\ge \P_{\delta_x}(V=0)\ge \exp\{-v(x)\}.
\end{equation}
It follows from \eqref{3.1} and \eqref{3.2} that
$$
\P_{\delta_x}(V=0)= \exp\{-v(x)\}\qquad \mbox{and  }\quad \P_{\delta_x}(V =\infty)=1-\exp\{-v(x)\}.
$$
By the dominated convergence theorem, we get that, for any $\theta>0$,
\begin{align*}
\lim_{t\to\infty}V_t(\theta v)(x)=-\log \P_{\delta_x}(e^{-\theta V})=-\log\P_{\delta_x}(V=0)=v(x).
\end{align*}
\hfill $\Box$

\section{Backward iterates}
It is clear that $V_t0=0$. It follows from \eqref{eq:v} that $V_tv=v$. Thus
$v$ and $0$ are two fixed points of $V_t$.

\begin{definition}
A family $(\eta_t,t\ge0)\subset \mathcal{B}^+(E)$
is called a family of backward iterates if $\eta_t\le v$ for all $t\ge 0$ and
$$
\eta_t(x)=V_s(\eta_{t+s})(x),\qquad t, s\ge0,\, x\in E.
$$
A  family  $(\eta_t, t\ge 0)$ of backward iterates
is said to be  non-trivial if, for some $t\ge0$, neither
$\eta_t=0, m$-a.e. nor $\eta_t=v, m$-a.e.
\end{definition}

It is well known that (see, for instance, \cite[1.44]{Rudin}),
there exists a metric $d$ on $C(E)$ such that $(C(E), d)$ is a complete metric space, and convergence in $(C(E), d)$
is equivalent to uniform convergence on each compact subset $K$ of $E$.

For any $a>0$, let ${\cal D}_a(E):=\{f\in {\cal B}^+(E): \|f\|_\infty\le a\}$.

\begin{lemma}\label{compact}
For any $t>0$, $V_t({\cal D}_a(E))$ is
a relatively compact subset of $C(E)$.
\end{lemma}

\noindent\textbf{Proof:}
Without loss of generality, we only prove the lemma for $t=1$ and $a\ge 1$.
We first show that, for any compact subset $K\subset E$,
 $\{V_1f: f\in {\cal D}_a(E)\}$ are equicontinuous on $K$.

Recall that
$$
V_1f(x)=T_1f(x)-\int_0^1 T_s(r(\cdot,V_{1-s}f))(x){\rm d}s.
$$
It is clear that for $f\in {\cal D}_a(E)$,
$$
|T_1f(x)-T_1f(y)|\le a\int_E|q(1, x,z)-q(1,y,z)|\,m({\rm d}z).
$$
Note that by \eqref{condition2}, we have
\begin{align}\label{e:rsnew2}
q(1, x,z)\le c_1\|\phi_0\|_\infty\psi_0(z).
\end{align}
Since $\psi_0\in L^1(E, m)$, for any $\epsilon>0$,  we can choose a compact set $\widetilde{K}\subset E$ such that
$$
2c_1\|\phi_0\|_\infty
\int_{\widetilde{K}^c}\psi_0(z)m({\rm d}z)<\frac{\epsilon}2.
$$
Using the continuity of $q(1, \cdot, \cdot)$ on $K\times \widetilde{K}$,
\eqref{e:rsnew2} and the fact that $\psi_0\in L^1(E, m)$,
we can find a
$\delta>0$ such that for any $x, y\in K$ with $|x-y|<\delta$,
$$
\int_{\widetilde{K}}|q(1, x,z)-q(1,y,z)|\,
m({\rm d}z)<\frac{\epsilon}2.
$$
Thus $\{T_1f: f\in {\cal D}_a(E)\}$
are equicontinuous on $K$.
By \eqref{2.6} and the fact that $V_tf(x)\le T_tf(x)\le e^{Mt}\|f\|_\infty$, we have that for any $f\in {\cal D}_a(E)$, $r(x,V_{1-s}f(x))\le 2M(e^{2M}a^2+e^{M}a).$
Thus
\begin{align*}
 &\left|\int_0^1 T_s(r(\cdot V_{1-s}f)(x){\rm d}s-\int_0^1 T_s(r(\cdot,V_{1-s}f)(y){\rm d}s\right|\\
 &\le
 [2M(e^{2M}a^2+e^{M}a)]
 \int_0^1 \int_E|q(s, x,z)-q(s, y,z)|\,
 m({\rm d}z){\rm d}s.
\end{align*}
Note that, for $\eta\in(0,1)$,
\begin{align*}
&\int_0^1 \int_E|q(s, x,z)-q(s, y,z)|\,
m({\rm d}z){\rm d}s\\
  &\le \int_0^\eta \int_Eq(s, x,z)+q(s, y,z)\,
   m({\rm d}z){\rm d}s+ \int_\eta^1 \int_E|q(s, x,z)-q(s, y,z)|\,
m({\rm d}z){\rm d}s\\
  &\le2e^M\eta+ \int_\eta^1 \int_E|q(s, x,z)-q(s, y,z)|\,
    m({\rm d}z){\rm d}s\\
  &\le2e^M\eta+ \int_\eta^1 \int_E\int_E|q(\frac{\eta}2, x,w)-q(\frac\eta2, y,w)|
  q(s-\frac\eta2, w, z)\,
   m({\rm d}w)\,m({\rm d}z){\rm d}s.
\end{align*}
For any $\epsilon>0$, we choose $\eta\in (0, 1)$ so that
$2e^M\eta<\frac\epsilon3$.
It follows from \eqref{density0} that, for any $s\in[\eta/2,1)$,
\begin{align}\label{e:rsnew}
q(s, x,z)\le (1+c(\eta/2))e^{\lambda_0}\|\phi_0\|_\infty\psi_0(z).
\end{align}
Hence
\begin{align*}
&\int_\eta^1 \int_E\int_E|q(\frac{\eta}2, x,w)-q(\frac\eta2, y,w)|
  q(s-\frac\eta2, w, z)\,
m({\rm d}w)\,m({\rm d}z){\rm d}s\\
&  \le (1+c(\eta/2))e^{\lambda_0}\|\phi_0\|_\infty\int_\eta^1
\int_E\int_E|q(\frac{\eta}2, x,w)-q(\frac\eta2, y,w)|
\,m({\rm d}w)\psi_0(z)\,m({\rm d}z){\rm d}s.
\end{align*}
Applying \eqref{e:rsnew} again, we can find a compact $\widetilde{K}\subset E$ such that
$$
(1+c(\eta/2))e^{\lambda_0}\|\phi_0\|_\infty\int_\eta^1
\int_E\int_{\widetilde{K}^c}|q(\frac{\eta}2, x,w)-q(\frac\eta2, y,w)|
\,m({\rm d}w)\psi_0(z)\,m({\rm d}z){\rm d}s<\frac\epsilon3.
$$
Using the continuity of  $q(\frac\eta2, \cdot, \cdot)$ on $K\times \widetilde{K}$,
\eqref{e:rsnew} and the fact that $\psi_0\in L^1(E, m)$, we can find a
$\delta>0$ such that for any $x, y\in K$ with $|x-y|<\delta$,
$$
(1+c(\eta/2))e^{\lambda_0}\|\phi_0\|_\infty\int_\eta^1
\int_E\int_{\widetilde{K}}|q(\frac{\eta}2, x,w)-q(\frac\eta2, y,w)|
\,m({\rm d}w)\psi_0(z)\,m({\rm d}z){\rm d}s<\frac\epsilon3.
$$
Thus $\{\int_0^1 T_s(r(\cdot,V_{1-s}f))(x){\rm d}s: f\in {\cal D}_a(E)\}$
are equicontinuous on $K$.
It follows that $\{V_1f: f\in {\cal D}_a(E) \}$
are equicontinuous on $K$.
In particular,
$V_1({\cal D}_a(E))\subset C(E)$.

Let $K_n$ be an increasing sequence of compact subsets of $E$.  Using the equicontinuity of $\{V_1f: f\in {\cal D}_a(E) \}$ on each $K_n$ and a routine diagonalization argument, we can easily show that any sequence of functions
in $V_1({\cal D}_a(E))$ contains a subsequence which converges in $C(E)$.
\hfill $\Box$

\begin{proposition}\label{lem:back}
There exists a non-trivial family of backward iterates.
\end{proposition}

\noindent\textbf{Proof:}
Let $\mathcal{D}=\{f\in\mathcal{B}^+(E):f\le  v\}$. For any $f\in{\cal D}$, $0\le V_sf\le V_sv\le v$. Thus $V_{t+s}f=V_{t}(V_sf)\in V_t({\cal D})$, which implies that $t\to V_t({\cal D})$ is decreasing.

For any $g, h\in {\cal D}$ and $t>0$, it is easy to see that
$\lambda\mapsto V_t(\lambda g+(1-\lambda)h)$, $\lambda\in [0, 1]$, is a continuous curve in  $V_t({\cal D})$ connecting $V_tg$ to $V_th$.
Thus $V_t({\cal D})$ is connected.
Hence $\mathcal{D}_\infty:=\bigcap_{n\in \mathbb{N}}V_n(\mathcal{D})$ is connected. Since $0,v\in \mathcal{D}_\infty$ and $v>0$, there exists an $\eta\in \mathcal{D}_\infty$, such that $0<\eta<v$ on a set of positive measure.
Thus, for every $n\in\mathbb{N}$, there exists $\eta_{n,n}\in \mathcal{D}$ such that $V_n(\eta_{n,n})=\eta$. Define $\eta_{n,j}=V_{n-j}(\eta_{n,n})$, $j=0,1,2,\cdots,n-1$. Note that $\eta_{n,0}=\eta$.

Since $0\le \eta_{n,j}\le v$ and $\eta_{n,j}=V_{1}(\eta_{n,j+1})$, it follows from
Lemma \ref{compact} that $\{\eta_{n,j}, n\ge 0, 0\le j<n\}$ is relatively compact
in $C(E)$. Thus
there exists a sequence $(\eta_{n_l,j})$ such that
 $\eta_j:=\lim_{l\to\infty}\eta_{n_l,j}\in C(E)\bigcap \mathcal{D}$ exists.
Since $\eta_{n_l,j}=V_{1}(\eta_{n_l,j+1})$, letting $l\to\infty$, we get that $\eta_j=V_1(\eta_{j+1})$.

Define $\eta_t:=V_{[t+1]-t}(\eta_{[t+1]})$, for $t\ge 0$. Then
$$
V_s(\eta_{t+s})=V_s V_{[t+s+1]-t-s}(\eta_{[t+s+1]})= V_{[t+s+1]-t}(\eta_{[t+s+1]})=V_{[t+1]-t}(\eta_{[t+1]})=\eta_t.
$$
It follows from  $\eta_{n,0}=\eta$ that $\eta_0=\eta$, which implies that
the family $\{\eta_t,t\ge 0\}$ of backward iterates is nontrivial.
\hfill $\Box$

\begin{lemma}\label{lem6}
If $0\le f\le v$, then
$$
V_tf(x)\ge e^{-at}T_tf(x),\quad t\ge 0, x\in E,
$$
where $a=2M(1+\|v\|_\infty)$.
\end{lemma}

\noindent\textbf{Proof:}
Using \eqref{2.6} with $H=1$ and the fact that $V_{t-s}f(x)\le V_{t-s}v(x)=v(x)\le \|v\|_\infty$
for any $0\le s\le t$ and any $x\in E$,
we get
$$
r(x,V_{t-s}f(x))\le 2M(1+\|v\|_\infty)V_{t-s}f(x)=aV_{t-s}f(x),\quad 0\le s\le t, x\in E.
$$
Recall that
$$
V_{t}f(x)=T_tf(x)-\int_0^tT_s(r(\cdot,V_{t-s}f(\cdot)))(x){\rm d}s.
$$
Thus we have
$$
V_{t}f(x)=e^{-at}T_tf(x)-\int_0^te^{-as}T_s(r(\cdot,V_{t-s}f(\cdot)))(x){\rm d}s+a\int_0^te^{-as}T_s (V_{t-s}f)(x){\rm d}s.
$$
Consequently we have
$$
V_{t}f(x)\ge e^{-at}T_tf(x),\quad t\ge 0, x\in E.
$$
\hfill $\Box$

\begin{lemma}\label{lem4}
If
$(\eta_t, t\ge 0)$
is a non-trivial family of backward iterates, then
$$
\lim_{t\to\infty}\|\eta_t\|_{\infty}=0.
$$
\end{lemma}

\noindent\textbf{Proof:}
Without loss of generality, we may and will assume that $m(\eta_0>0)>0$ and
$m(\eta_0<v)>0$. We claim that $\la \eta_t,\psi_0\ra_m\to 0$ as $t\to\infty$.
Otherwise,
there exist a sequence $t_j\uparrow\infty$ and a constant $c_0>0$ such that
 $$
 \la \eta_{t_j},\psi_0\ra_m>c_0.
 $$
Fix $s>1$ large enough such that
$1-ce^{-\gamma s}>0$,
where $c$ and $\gamma$ are the constants in \eqref{T_t} with $\delta=1$.
Then, for $j$ large enough so that $t_j>s$, we have
 \begin{equation}\label{3.6}
   \eta_0(x)=V_{t_j}(\eta_{t_j})(x)=V_{t_j-s}(V_s(\eta_{t_j}))(x),\quad x\in E.
 \end{equation}
By Lemma \ref{lem6}, we have
$$
V_s(\eta_{t_j})(x)\ge e^{-as}T_s(\eta_{t_j})(x)\ge e^{-as}(1-ce^{-\gamma s})e^{\lambda_0s}\la\eta_{t_j},\psi_0 \ra_m\phi_0(x)\ge c_0e^{-as}(1-ce^{-\gamma s})e^{\lambda_0s} \phi_0(x).
$$
It follows from \eqref{de:R1} and \eqref{1.10} that
\begin{equation}\label{3.3}
  v(x)=V_1v(x)\le T_1v(x)\le (1+c)e^{\lambda_0}\la v,\psi_0 \ra_m\phi_0(x).
\end{equation}
Thus
\begin{equation}\label{3.5}
  V_s(\eta_{t_j})(x)\ge c_0e^{-as}(1-ce^{-\gamma s})e^{\lambda_0s} (1+c)^{-1}e^{-\lambda_0}\la v,\psi_0 \ra_m^{-1} v(x)=:C_sv(x).
\end{equation}
Combining \eqref{3.5}, \eqref{3.6} and \eqref{3.11}, we get
$$
\eta_0(x)\ge \lim_{j\to\infty}V_{t_j-s}(C_sv)(x)=v(x),
$$
which contradicts  the definition of $\eta_0$.
Thus the claim above is true.

Note that, for any $s\ge1$ and $t>0$,
\begin{equation}\label{3.8}
  \eta_t(x)=V_s(\eta_{t+s})(x)\le T_s (\eta_{t+s})(x)
  \le (1+ce^{-\gamma s})e^{\lambda_0s}\la \eta_{t+s},\psi_0\ra_m\phi_0(x).
\end{equation}
Thus
$$
\|\eta_t\|_\infty\le (1+c)e^{\lambda_0s}\la \eta_{t+s},\psi_0\ra_m\|\phi_0\|_\infty\to 0,
$$
as $t\to\infty$.
The proof is now complete.
\hfill $\Box$

\smallskip

\begin{lemma}\label{lem8}
If $(\eta_t, t\ge 0)$
is a non-trivial family of backward iterates, then
there exist  $\{h_t: t\ge 0\}\subset \mathcal{B}_b(E)$ such that
$$
\eta_t(x)=(1+h_t(x))\la\eta_t,\psi_0\ra_m\phi_0(x),\quad t\ge 0, x\in E,
$$
and
$$
\lim_{t\to\infty}\|h_t\|_\infty=0.
$$
Moreover,
\begin{equation}\label{3.10}
  \frac{\la\eta_t,\psi_0\ra_m}{\la\eta_{t+s},\psi_0\ra_m}\ge e^{\lambda_0s}, \,\forall s, t\ge 0
  \quad\mbox{and} \quad
  \lim_{t\to\infty}\frac{\la\eta_t,\psi_0\ra_m}{\la\eta_{t+s},\psi_0\ra_m}=e^{\lambda_0s},\,\forall s\ge 0.
\end{equation}
\end{lemma}

\noindent\textbf{Proof:}
Since $\la\eta_t,\psi_0\ra_m\phi_0(x)>0$, we can define
$$
h_t(x):=\frac{\eta_t(x)}{\la\eta_t,\psi_0\ra_m\phi_0(x)}-1, \quad t\ge 0, x\in E.
$$
It follows from \eqref{T_t} and Lemma \ref{lem1.2}  that for $s>1$,
\begin{align}
  \eta_t(x)&=V_s(\eta_{t+s})(x)=T_s(\eta_{t+s})(x)-R_{\eta_{t+s}}(s,x)\label{5.1}\\
  &\ge e^{\lambda_0s}(1-ce^{-\gamma s}
  -\|g_{\eta_{t+s}}(s)\|_\infty)\la\eta_{t+s},\psi_0\ra_m\phi_0(x)\label{3.7},
\end{align}
where $c$ and $\gamma$ are the constants in \eqref{T_t} with $\delta=1$.
In the remainder of this proof, we assume that $s$ is large enough such that  $1-ce^{-\gamma s}>0$.
By Lemma \ref{lem4} and \eqref{1.22}, $\lim_{t\to\infty}\|g_{\eta_{t+s}}(s)\|_\infty=0$.
Thus, for large enough $t$,
$1-ce^{-\gamma s}-\|g_{\eta_{t+s}}(s)\|_\infty>0$.
It follows from \eqref{3.8} and \eqref{3.7} that
\begin{align*}
 -\frac{2ce^{-\gamma s}+\|g_{\eta_{t+s}}(s)\|_\infty}{1+ce^{-\gamma s}}
 &\le h_t(x)=\frac{\eta_t(x)}{\la\eta_t,\psi_0\ra_m\phi_0(x)}-1
   \le\frac{2ce^{-\gamma s}+\|g_{\eta_{t+s}}(s)\|_\infty}{1-ce^{-\gamma s}-\|g_{\eta_{t+s}}(s)\|_\infty} .
\end{align*}
Letting $t\to\infty$ and then $s\to\infty$, we get that
$$
\lim_{t\to\infty}\|h_t\|_\infty=0.
$$
It follows from \eqref{5.1} that
\begin{align*}
  \la\eta_t,\psi_0\ra_m&=\la T_s(\eta_{t+s}),\psi_0\ra_m-\la R_{\eta_{t+s}}(s,x),\psi_0\ra_m\\
  &=e^{\lambda_0s}\la\eta_{t+s},\psi_0\ra_m-e^{\lambda_0s}\la\eta_{t+s},\psi_0\ra_m\la g_{\eta_{t+s}}(s)\phi_0,\psi_0\ra_m,
  \end{align*}
which implies that
$$
0\le e^{\lambda_0s}-\frac{\la\eta_t,\psi_0\ra_m}{\la\eta_{t+s},\psi_0\ra_m}
=e^{\lambda_0s}\la g_{\eta_{t+s}}(s)\phi_0,\psi_0\ra_m\le e^{\lambda_0s}\|g_{\eta_{t+s}}(s)\|_\infty
\to0,
$$
as $t\to\infty$, where the last limit follows from \eqref{1.22} and Lemma \ref{lem4}.
The proof is now complete.
\hfill $\Box$

\begin{remark}\label{Rem-4.7}
Let $L(t):=e^{\lambda_0t}\la\eta_t,\psi_0\ra_m$.
It follows from \eqref{3.10} that $t\to L(t)$ is nondecreasing and, for any $s\ge0$,
$$
\lim_{t\to\infty}\frac{L(t+s)}{L(t)}=1.
$$
Therefore, $l_0:=\lim_{t\to\infty}L(t)\in(0,\infty]$ exists.
\end{remark}

\section{Seneta-Heyde norming for $M$}

\begin{lemma}\label{lem10}
If $f\in {\cal B}^+(E)$ satisfies $f(x)\le v(x)$ and $f(x)=V_tf(x)$ for all $t>0$ and $x\in E$,
then either $f(x)=v(x)$ for all $x\in E$ or $f(x)=0$ for all $x\in E$.
\end{lemma}

\noindent\textbf{Proof:}
Let $s_0>1$ be such that
$1-ce^{-\gamma s_0}>0$, where $c$ and $\gamma$ are the constants from \eqref{T_t} with $\delta=1$.
By Lemma \ref{lem6} and \eqref{3.3}, we have
$$
f(x)=V_{s_0}f(x)\ge e^{-as_0}T_{s_0}f(x)
\ge e^{-as_0}(1-ce^{-\gamma s_0})e^{\lambda_0s_0}\la f,\psi_0\ra_m\phi_0(x)
\ge c_0\la f,\psi_0\ra_mv(x)
$$
for some constant $c_0>0$.
If $\la f,\psi_0\ra_m>0$, then by \eqref{3.11} we have
$$
f(x)=\lim_{t\to\infty}V_tf(x)\ge \lim_{t\to\infty}V_t(
c_0\la f,\psi_0\ra_mv)(x)=v(x),
$$
which implies that $f(x)=v(x)$ for all $x\in E$.
If $\la f,\psi_0\ra_m=0$, then
$$
f(x)\le T_1f(x)\le (1+c)e^{\lambda_0}\la f,\psi_0\ra_m \phi_0(x)=0,
$$
which implies that $f(x)=0$ for all $x\in E$.
The proof is now complete.
\hfill $\Box$

\begin{theorem}\label{Thm:main1}
Let $(\eta_t, t\ge 0)$
be a non-trivial family of backward iterates and
define $\gamma_t:=\langle\eta_t,\psi_0\rangle_m$.
Then there is a non-degenerate random variable $W$ such that for
any nonzero $\mu\in\mathcal{M}_F(E)$,
$$
\lim_{t\to\infty}\gamma_t\langle \phi_0,X_t\rangle=W,\qquad
\mbox{a.s.-}\P_{\mu}
$$
and
\begin{equation}\label{3.9}
\P_{\mu}(W=0)=\exp\{-\la v,\mu\ra\},\qquad \P_{\mu}(W<\infty)=1.
\end{equation}
\end{theorem}

\noindent\textbf{Proof:}
By the definition of $\eta_t$, for any nonzero $\mu\in\mathcal{M}_F(E)$,
$$
\P_{\mu}(\exp\{-\eta_{t+s},X_{t+s}\}|\mathcal{G}_t)=\P_{X_t}(\exp\{-\eta_{t+s},X_{s}\})=\exp\{-\langle \eta_t,X_t\rangle \},
$$
which implies that $\{\exp\{-\langle \eta_t,X_t\rangle \},t\ge 0\}$ is a
nonnegative martingale.
Thus, by the martingale convergence theorem,
$\lim_{t\to\infty}\exp\{-\langle \eta_t,X_t\rangle \}$ exists $\P_{\mu}$ almost surely
and hence
$W:=\lim_{t\to\infty}\langle \eta_t,X_t\rangle\in[0,\infty]$ exists $\P_{\mu}$
almost surely.

It follows from Lemma \ref{lem8} that
$$
(1-\|h_t\|_\infty)\gamma_t\langle \phi_0,X_t\rangle
\le\la\eta_t,X_t\ra\le (1+\|h_t\|_\infty)\gamma_t\langle \phi_0,X_t\rangle.
$$
Since $\lim_{t\to\infty}\|h_t\|_\infty=0$, we have $1-\|h_t\|_\infty>0$ for $t$ large enough.
Thus for large $t$,
$$
(1+\|h_t\|_\infty)^{-1}\la\eta_t,X_t\ra
\le\gamma_t\langle \phi_0,X_t\rangle\le (1-\|h_t\|_\infty)^{-1}\la\eta_t,X_t\ra.
$$
Letting $t\to\infty$, we get that
$$
\lim_{t\to\infty}\gamma_t\langle \phi_0,X_t\rangle=W,\qquad
\mbox{a.s.-}\P_{\mu}.
$$

Define \begin{equation}\label{def-Phi}\Phi(s,x):=-\log \P_{\delta_x}\exp\{-sW\}.\end{equation}
Then
\begin{align*}
&-\log\P_{\mu}\exp\{-sW\}=\lim_{t\to\infty}-\log\P_{\mu}\exp\{-s\gamma_t\langle \phi_0,X_t\rangle\}\\
&=\lim_{t\to\infty}\la-\log\P_{\delta_x}\exp\{-s\gamma_t\langle \phi_0,X_t\rangle\},\mu\ra
=\la \Phi(s,\cdot),\mu\ra.
\end{align*}
Put $\Phi_\infty(x):=\lim_{s\to\infty}\Phi(s,x)$ and $\Phi_0(x)=\lim_{s\to0}\Phi(s,x)$.
Then
$$
\P_{\mu}(W=0)=\lim_{s\to\infty}\P_{\mu}\exp\{-sW\}=\exp\{-\la \Phi_\infty,\mu\ra\}
$$
and
$$
\P_{\mu}(W<\infty)=\lim_{s\to0}\P_{\mu}\exp\{-sW\}=\exp\{-\la\Phi_0,\mu\ra\}.
$$

For any $s,t>0$, we have
\begin{align*}
   &\exp\{-\Phi(s,x)\} = \lim_{u\to\infty}\P_{\delta_x}(\exp\{-s\gamma_{t+u}\la\phi_0,X_{t+u}\ra\}) \\
   &= \lim_{u\to\infty}\P_{\delta_x}\P_{X_t}(\exp\{-s\gamma_{t+u}\gamma_u^{-1}\gamma_u\la\phi_0,X_{u}\ra\})\\
   &=\P_{\delta_x}\P_{X_t}\exp\{-se^{-\lambda_0t}W\}=\P_{\delta_x}\exp\{-\la \Phi(se^{-\lambda_0t},\cdot),X_t\ra\},
\end{align*}
which implies that
\begin{equation}\label{Phi}
\Phi(s,x)=V_t(\Phi(se^{-\lambda_0t}, \cdot))(x).
\end{equation}
Thus, letting $s\to\infty$ and $s\to0$, we get respectively
$$
\Phi_\infty(x)=V_t(\Phi_\infty)(x),\qquad \Phi_0(x)=V_t(\Phi_0)(x).
$$
Since $s\to \Phi(s,x)$ is increasing, we have $\Phi_0(x)\le \Phi(1,x)\le \Phi_\infty(x)$.
Note that
$$
\Phi(1,x)=-\log\P_{\delta_x}\exp\{-W\}
=\lim_{t\to\infty}-\log\P_{\delta_x}\exp\{-\langle \eta_t,X_t\rangle \}
=\eta_0(x).
$$
On the other hand,
$$
\P_{\delta_x}(W=0)\ge \P_{\delta_x}(\exists t>0, \|X_t\|=0)=\exp\{-v(x)\},
$$
which implies that, for all $x\in E$,
$
\Phi_{\infty}(x)\le v(x).
$
Thus, $\Phi_0(x)\le \eta_0(x)\le \Phi_\infty(x)\le v(x)$.
It follows from Lemma \ref{lem10} that
$$
\Phi_0=0, \quad \Phi_\infty=v.
$$
The proof is now complete.
\hfill $\Box$

Now Theorem \ref{main1} follows immediately from Lemma \ref{lem8} and Theorem \ref{Thm:main1}.

\begin{proposition}
Let $(\eta_t, t\ge 0)$
be a non-trivial family of backward iterates and $W$ be the limit in Theorem \ref{Thm:main1} corresponding to $(\eta_t, t\ge 0)$.
Then a family $(\eta^*_t, t\ge 0)$
 is a non-trivial family of backward iterates  if and only if there exists a positive number $a$ such that
$$\eta^*_t(x)=\Phi(ae^{-\lambda_0t},x),\quad t\ge 0, x\in E,$$
where $\Phi$ is defined in \eqref{def-Phi}.
In particular,
\begin{equation}\label{eq:eta}
  \eta_t(x)=\Phi(e^{-\lambda_0t},x),\quad t\ge 0, x\in E.
\end{equation}
\end{proposition}

\noindent\textbf{Proof:}
For any $a>0$, by \eqref{Phi}, we have
$$
V_s(\Phi(ae^{-\lambda_0(t+s)}, \cdot))(x)=\Phi(ae^{-\lambda_0t},x).
$$
Thus $(\Phi(ae^{-\lambda_0t},x), t\ge0)$ is a non-trivial family of backward iterates.

If $(\eta^*_t, t\ge 0)$ is a non-trivial family of backward iterates,
then it follows from Lemma \ref{lem8} that, for any $s\ge 0$, $(\eta^*_{t+s},t\ge0)$ is also a  non-trivial family of backward iterates.
Let $W^{(s)}$ be the limit in Theorem \ref{Thm:main1} corresponding to $(\eta^*_{t+s},t\ge0)$.
By \eqref{3.9}, we get that, for any $s\ge0$,
$$
\{W>0 \}=\{W^{(s)}>0\}=\{\forall\, t>0, \,\|X_t\|>0\},\qquad \mbox{a.s.-}\P_{\mu}.
$$
Thus, we have that, for any $\omega\in \{W>0\}\cap \{W^{(s)}>0\}$,
$$
\frac{W^{(s)}(\omega)}{W(\omega)}
=\lim_{t\to\infty}\frac{\la \eta^*_{t+s},\psi_0\ra_m}{\la \eta_t,\psi_0\ra_m}
:=e(s)>0,
$$
where $e(s)$ is deterministic. Therefore,
\begin{equation}\label{4.21}
  W^{(s)}=e(s)W, \qquad \mbox{a.s.-}\P_{\mu}.
\end{equation}
By \eqref{3.10}, we have that, for any $s,r\ge 0$,
$$
W^{(s)}=\lim_{t\to\infty}\frac{\la \eta^*_{t+s},\psi_0\ra_m}{\la \eta^*_{t+r},\psi_0\ra_m}W^{(r)}
=e^{-\lambda_0(s-r)}W^{(r)},\qquad \mbox{a.s.-}\P_{\mu}.
$$
It follows that
$e(s)=e^{-\lambda_0(s-r)}e(r)$, which implies that there
exists a constant $a>0$ such that
$$e(s)=ae^{-\lambda_0s}.$$
Note that, for any $s\ge 0$ and $x\in E$,
\begin{align*}
  \eta^*_s(x) &= \lim_{t\to\infty}-\log \P_{\delta_x}\exp\{-\la\eta^*_{t+s},X_t\ra\}
  =-\log \P_{\delta_x}\exp\{-W^{(s)}\} \\
  &=  -\log \P_{\delta_x}\exp\{-ae^{-\lambda_0s}W\}=\Phi(ae^{-\lambda_0s},x).
\end{align*}
The proof is now complete.
\hfill $\Box$

Recall from Remark \ref{Rem-4.7} that
$l_0=\lim_{t\to\infty}e^{\lambda_0t}\gamma_t\in(0,\infty]$.

\begin{proposition}\label{prop1}
Let $(\eta_t, t\ge 0)$
be a non-trivial family of backward iterates and $W$ be the limit in Theorem \ref{Thm:main1} corresponding to $(\eta_t, t\ge 0)$.
\begin{enumerate}
  \item[(1)]  If $l_0<\infty$, then $\P_{\mu}W<\infty$ for any $\mu\in \mathcal{M}_F(E)$.
  \item[(2)] If  $\P_{\mu}W<\infty$ for some nonzero $\mu\in \mathcal{M}_F(E)$, then $l_0<\infty$.
\end{enumerate}
Moreover,
\begin{equation}\label{exp:W}
  \P_{\mu}W=\lim_{t\to\infty}\gamma_t\P_{\mu}\la\phi_0,X_t\ra=l_0\la\phi_0,\mu\ra.
\end{equation}
\end{proposition}

\noindent\textbf{Proof:}
(1) Since $W=\lim_{t\to\infty}\gamma_t\la\phi_0,X_t\ra$,
we have by Fatou's lemma that
$$\P_{\mu}W\le \lim_{t\to\infty}\P_{\mu}\gamma_t\la\phi_0,X_t\ra=\lim_{t\to\infty}\gamma_te^{\lambda_0t}\la\phi_0,\mu\ra=l_0\la\phi_0,\mu\ra.$$
Thus $l_0<\infty$ implies $\P_{\mu}W<\infty$ for any $\mu\in{\cal M}_F(E)$.

(2) It follows from \eqref{eq:eta} and Lemma \ref{lem8} that
\begin{align}\label{5.5}
  \lim_{s\to0}\frac{\la\Phi(s,\cdot),\mu\ra}{s}
   &= \lim_{t\to\infty} e^{\lambda_0t}\la\Phi(e^{-\lambda_0t},\cdot),\mu\ra=\lim_{t\to\infty}e^{\lambda_0t}\la\eta_t,\mu\ra\nonumber\\
   &=\lim_{t\to\infty}e^{\lambda_0t}\gamma_t\la(1+h_t)\phi_0,\mu\ra=l_0\la\phi_0,\mu\ra.
\end{align}

If $\P_{\mu}W<\infty$ for some nonzero $\mu\in{\cal M}_F(E)$, then
\begin{eqnarray*}
  \P_{\mu}W=\lim_{s\to0}\frac{1-\P_{\mu}(\exp\{-sW\})}{s}=\lim_{s\to0}\frac{\la\Phi(s,\cdot),\mu\ra}{s}=l_0\la\phi_0,\mu\ra,
\end{eqnarray*}
which implies $l_0<\infty$.
The proof is now complete.
\hfill $\Box$

It was shown in \cite[Section 1.10, Lemma 1]{AN72} that,
if $Y$ is a non-negative random variable with $EY<\infty$, then
$$
EY\log^+Y<\infty
$$ if and only if, for some $a_0>0$
$$
\int_0^{a_0} s^{-2}E(e^{-s Y}-1+s Y)\,{\rm d}s<\infty.
$$

Recall that
$l(x):=\int_1^\infty \theta\log \theta \,n^{\phi_0}(x, {\rm d}\theta).$
\begin{lemma}\label{lem:llogl}
$\int_E\psi_0(x)l(x)m({\rm d}x)<\infty$
if and only if, for any $t>0$,
$$ \int_E\psi_0(x)\P_{\delta_x}\Big[\la\phi_0,X_t\ra\log^+\la\phi_0,X_t\ra\Big]m({\rm d}x)<\infty.$$
\end{lemma}

\noindent\textbf{Proof:}
Without loss of generality, we only prove the result for $t=1$.
Note that
$$
\int_E\psi_0(x)\P_{\delta_x}\Big[\la\phi_0,X_1\ra\log^+\la\phi_0,X_1\ra\Big]m({\rm d}x)<\infty
$$
if and only if, for some $a_0>0$
\begin{equation}\label{6.1}
  \int_E\psi_0(x)m({\rm d}x)\int_0^{a_0} s^{-2}\P_{\delta_x}(e^{-s\la\phi_0,X_1\ra}-1+s\la\phi_0,X_1\ra)\,{\rm d}s<\infty.
\end{equation}
Put $R(f)(x):=R_f(1,x)=T_1f(x)-V_1f(x)$.  Then, we have
\begin{align*}
  &\P_{\delta_x}(e^{-s\la\phi_0,X_1\ra}-1+s\la\phi_0,X_1\ra)\\
  & = \exp\{-V_1(s\phi_0)(x)\}-1+ V_1(s\phi_0)(x)+R(s\phi_0)(x)\\
  &\le   \frac{1}{2}V_1(s\phi_0)(x)^2+R(s\phi_0)(x)\\
  &\le \frac{1}{2}T_1(s\phi_0)(x)^2+R(s\phi_0)(x)\\
  &= \frac{1}{2}e^{2\lambda_0}s^2\phi_0(x)^2+R(s\phi_0)(x).
\end{align*}
On the other hand,
$$
\P_{\delta_x}(e^{-s\la\phi_0,X_1\ra}-1+s\la\phi_0,X_1\ra) = \exp\{-V_1(s\phi_0)(x)\}-1+T_1(s\phi_0)(x)\ge
R(s\phi_0)(x).
$$
Thus,
\eqref{6.1} holds if and only if, for some $a_0>0$,
\begin{equation}\label{6.2}
  \int_0^{a_0} s^{-2}\la R(s\phi_0),\psi_0\ra_m{\rm d}s<\infty.
\end{equation}
By \eqref{2.2}, we have
\begin{align*}
&  \la r(\cdot, V_{1-t}(s\phi_0)),\psi_0\ra_m\\
&= \la \beta V_{1-t}(s\phi_0)^2,\psi_0\ra_m
  +\Big\la \Big(\int_0^{\phi_0^{-1}}+\int_{\phi_0^{-1}}^\infty\Big)
  \Big(e^{-\theta V_{1-t}(s\phi_0)}-1+\theta V_{1-t}(s\phi_0)\Big)n(\cdot,{\rm d}\theta),
  \psi_0\Big\ra_m\\
  &=: J_1(s,t)+J_2(s,t)+J_3(s,t).
\end{align*}
Thus, by \eqref{IF}, we have
\begin{align*}
&  \int_0^{a_0} s^{-2}\la R(s\phi_0),\psi_0\ra_m {\rm d}s\\
&=\int_0^{a_0} s^{-2}{\rm d}s\int_0^1 e^{\lambda_0t}\la r(\cdot, V_{1-t}(s\phi_0)),\psi_0\ra_m\,{\rm d}t\\
  &=\int_0^{a_0} s^{-2}{\rm d}s\int_0^1 e^{\lambda_0t} J_1(s,t)\,{\rm d}t+\int_0^{a_0} s^{-2}{\rm d}s\int_0^1 e^{\lambda_0t} J_2(s,t)\,{\rm d}t\\
&\qquad+\int_0^{a_0} s^{-2}{\rm d}s\int_0^1 e^{\lambda_0t} J_3(s,t)\,{\rm d}t.
\end{align*}
Since
\begin{equation}\label{6.4}
  V_{1-t}(s\phi_0)(x)\le T_{1-t}(s\phi_0)(x)\le s e^{\lambda_0}\phi_0(x),
\end{equation}
we have that $J_1(s,t)\le Me^{2\lambda_0}s^2\|\phi_0\|_\infty$. Thus
\begin{equation}\label{J1}
  \int_0^{a_0} s^{-2}{\rm d}s\int_0^1 e^{\lambda_0t}J_1(s,t){\rm d}t<\infty.
\end{equation}
Note that
$$
J_2(s,t)\le
\Big\la\int_0^{\phi_0^{-1}}
\theta^2V_{1-t}(s\phi_0)^2n(\cdot,{\rm d}\theta),\psi_0\Big\ra_m
\le e^{2\lambda_0}s^2\Big\la\int_0^{1}
\theta^2n^{\phi_0}(\cdot,{\rm d}\theta),\psi_0\Big\ra_m
\le Cs^2,$$
which implies that
\begin{equation}\label{J2}
  \int_0^{a_0} s^{-2}{\rm d}s\int_0^1 e^{\lambda_0t}J_2(s,t){\rm d}t<\infty.
\end{equation}
Now we deal with $J_3$.
By Lemma \ref{lem6} and \eqref{1.6} with $s_0>1$ large enough
such that $c(1)e^{-\gamma(1)s_0}<\frac12$,
$$
v(x)=V_{s_0}v(x)\ge e^{-as_0}T_{s_0}v(x)\ge
\frac12e^{(\lambda_0-a)s_0}\phi_0(x).
$$
We put $A_0=\frac12e^{(\lambda_0-a)s_0}$.
Hence by Lemma \ref{lem6}, for any $s\le A_0$,
\begin{equation}\label{6.5}
  V_{1-t}(s\phi_0)(x)\ge e^{-a}T_{1-t}(s\phi_0)(x)\ge e^{-a}s\phi_0(x).
\end{equation}
Thus, combining \eqref{6.4} and \eqref{6.5}, there exist $C_1,C_2>0$ such that for any $s\le A_0$,
\begin{align*}
& \Big\la\int_{\phi_0^{-1}}^\infty \Big(e^{-C_1\theta s\phi_0}-1+C_1 \theta s\phi_0\Big)n(\cdot,{\rm d}\theta),\psi_0\Big\ra_m
\le J_3(s,t)\\
  &\le \Big\la\int_{\phi_0^{-1}}^\infty \Big(e^{-C_2 \theta s\phi_0}-1+C_2\theta s\phi_0\Big)n(\cdot,{\rm d}\theta ),\psi_0\Big\ra_m.
\end{align*}
Note that for any $C>0$, we have
\begin{align*}
  &
  \int_0^{A_0}s^{-2}{\rm d}s\int_0^1 e^{\lambda_0t}
  \Big\la\int_{\phi_0^{-1}}^\infty
\Big(e^{-C\theta s\phi_0}-1+C\theta s\phi_0\Big)n(\cdot,{\rm d}\theta),\psi_0\Big\ra_m {\rm d}t \\
  &= \int_0^1 e^{\lambda_0t} {\rm d}t\int_E \psi_0(x)m({\rm d}x)
 \int_{1}^\infty n^{\phi_0}(x,{\rm d}\theta)
  \int_0^{A_0}
s^{-2}\Big(e^{-C\theta s}-1+C\theta s\Big){\rm d}s\\
  &=\int_0^1 e^{\lambda_0t} {\rm d}t\int_E \psi_0(x)m({\rm d}x)\int_{1}^\infty
\theta n^{\phi_0}(x,{\rm d}\theta) \int_0^{A_0\theta}
  s^{-2}\Big(e^{-Cs}-1+Cs \Big){\rm d}s.
\end{align*}
Since
$$
 \lim_{\theta\to\infty}\frac{
\int_0^{A_0\theta }
s^{-2}\Big(e^{-Cs}-1+Cs\Big){\rm d}s}{\log \theta}=C,
$$
we have
$$
\int_0^{A_0}
s^{-2}{\rm d}s\int_0^1 e^{\lambda_0t}J_3(s,t){\rm d}t<\infty\Longleftrightarrow\int_E l(x)\psi_0(x)m({\rm d}x)<\infty.
$$
Now the conclusion follows immediately.
\hfill $\Box$

Recall that
$$
M_\infty=\lim_{t\to\infty}M_t=\lim_{t\to\infty}e^{-\lambda_0t}\la\phi_0,X_t\ra.
$$
\begin{proposition}\label{prop2}
If $\int_E\psi_0(x)l(x)m({\rm d}x)<\infty$,
then for any nonzero $\mu\in{\cal M}_F(E)$, $M_\infty$ is non-degenerate under $\P_{\mu}$
 and $\P_{\mu}M_\infty=\la\phi_0,\mu\ra$.
\end{proposition}

\noindent\textbf{Proof:}
Suppose $\mu\in{\cal M}_F(E)$ is nonzero and fixed. For any $\theta>0$, put
$$\Psi_t(\theta,x):=-\log\P_{\delta_x}(\exp\{-\theta M_t\})\qquad\mbox{and}\quad \Psi(\theta,x):=-\log\P_{\delta_x}(\exp\{-\theta M_\infty\}).$$
Then for any $x\in E$, $\Psi_t(\theta,x)$ is non-increasing in $t$.
By the dominated convergence theorem and monotone convergence theorem, we have
\begin{equation}\label{M-infinity}-\log\P_{\mu}(\exp\{-\theta M_\infty\})=\lim_{n\to \infty}\la\Psi_n(\theta,\cdot),\mu\ra=\la\Psi(\theta,\cdot),\mu\ra.\end{equation}
We claim that there exists some $\theta>0$ such that
\begin{equation}\label{5.13}
  \la\Psi(\theta,\cdot),\psi_0\ra_m>0.
\end{equation}
By the Markov property of $X$, we have
\begin{align}\label{5.14}
  &\Psi_{t+s}(\theta,x)=-\log\P_{\delta_x}\P_{X_t}(\exp\{-\theta e^{-\lambda_0t}M_s\})\nonumber\\
  &=-\log\P_{\delta_x}\exp\{-\la\Psi_s(\theta e^{-\lambda_0t}),X_t\ra\}=V_t(\Psi_s(\theta e^{-\lambda_0t}))(x).
\end{align}
Letting $s\to\infty$, we get
$$
\Psi(\theta e^{\lambda_0t},x)=V_t(\Psi(\theta))(x).
$$
If \eqref{5.13} holds, then by  \eqref{1.6}, we get that, for $t>1$ large enough,
$$
\P_{\delta_x}\la\Psi(\theta),X_t\ra=T_t\Psi(\theta)(x)>0,
$$
which implies that
$\P_{\delta_x}(\la\Psi(\theta),X_t\ra>0)>0$. Hence we have that, for any $x\in E$,
 $\Psi(\theta e^{\lambda_0t},x)=V_t(\Psi(\theta))(x)>0$.
Thus, by \eqref{M-infinity},
 $$\P_{\mu} \exp\{-\theta e^{\lambda_0t}M_\infty\}=\exp\{-\la \Psi(\theta e^{\lambda_0t}),\mu\ra\}<1,$$
which implies that
$\P_{\mu}(M_\infty=0)<1$.

Now we prove claim \eqref{5.13}.
Put $R(f)(x):=R_f(1,x).$
It follows from \eqref{5.14} and \eqref{de:R1} that
\begin{align}
& \la \Psi_{n}(\theta,\cdot),\psi_0\ra_m  = \la V_1(\Psi_{n-1}(\theta e^{-\lambda_0})),\psi_0\ra_m \nonumber\\
   &=  \la T_1(\Psi_{n-1}(\theta e^{-\lambda_0})),\psi_0\ra_m-\la R(\Psi_{n-1}(\theta e^{-\lambda_0})),\psi_0\ra_m\nonumber\\
   &= e^{\lambda_0}\la \Psi_{n-1}(\theta e^{-\lambda_0}),\psi_0\ra_m-\la R(\Psi_{n-1}(\theta e^{-\lambda_0})),\psi_0\ra_m\nonumber\\
   &=e^{\lambda_0(n-1)}\la \Psi_{1}(\theta e^{-(n-1)\lambda_0}),\psi_0\ra_m-\sum_{k=1}^{n-1} e^{\lambda_0(k-1)}\la R(\Psi_{n-k}(\theta e^{-\lambda_0k})),\psi_0\ra_m.
\end{align}
Note that, by Jensen's inequality, we have
$$\Psi_{n-k}(\theta e^{-\lambda_0k})(x)=-\log\P_{\delta_x}\exp\{-\theta e^{-\lambda_0k}M_{n-k}\}\le \P_{\delta_x}\theta e^{-\lambda_0k}M_{n-k}=\theta e^{-\lambda_0k}\phi_0(x).$$
By the dominated convergence theorem, we get
$$\lim_{n\to\infty}e^{\lambda_0(n-1)}\la \Psi_{1}(\theta e^{-(n-1)\lambda_0}),\psi_0\ra_m=\theta\la \P_{\delta_\cdot}M_1,\psi_0\ra_m=\theta.$$
Thus we have
\begin{equation}\label{5.15}
  \la\Psi(\theta,\cdot),\psi_0\ra_m
  =\lim_{n\to\infty}\la \Psi_{n}(\theta,\cdot),\psi_0\ra_m
  \ge \theta-\lim_{n\to\infty}\sum_{k=1}^{n-1} e^{\lambda_0(k-1)}\la R(\theta e^{-\lambda_0k}\phi_0),\psi_0\ra_m.
\end{equation}
Since $t\to \la R(\theta e^{-\lambda_0t}\phi_0),\psi_0\ra_m$ is decreasing, we have
\begin{align*}
  &\sum_{k=1}^{\infty} e^{\lambda_0(k-1)}\la R(\theta
  e^{-\lambda_0k}\phi_0),\psi_0\ra_m\\
  &\le \int_0^\infty e^{\lambda_0t}\la R(\theta e^{-\lambda_0t}\phi_0),\psi_0
 \ra_m\,{\rm d}t \\
   &= \lambda_0^{-1}\int_0^1 s^{-2}\la R(s\theta \phi_0),\psi_0
  \ra_m\,{\rm d}s\\
  &\le \lambda_0^{-1}\int_0^1 s^{-2}\la T_1(s\theta \phi_0)-
 (1-\exp\{-V_1(s\theta \phi_0)\}),\psi_0\ra_m\,{\rm d}s\\
  &=\lambda_0^{-1}\int_E\psi_0(x)m({\rm d}x)\P_{\delta_x}\int_0^1 s^{-2}\Big(s\theta \la\phi_0,X_1\ra-1+\exp\{-s\theta \la\phi_0,X_1\ra\}\Big)\,{\rm d}s.
\end{align*}
Since $e^{-s}-1+s\le s\wedge(s^2/2)$,
there exists $C>0$ such that for any $r\ge 0$,
\begin{align*}
  &\int_0^1 s^{-2}\Big(rs-1+\exp\{-rs\}\Big)\,{\rm d}s
   = r\int_0^{r}s^{-2}\Big(s-1+\exp\{-s\}\Big)\,{\rm d}s\\
   &\le
   \frac{1}{2}r^2{\rm I}_{r\le 2}+r\left(
 1+\int_2^r s^{-1}\,{\rm d}s\right){\rm I}_{r>2}
 \le \frac{1}{2}r^2{\rm I}_{r\le 2}+Cr\log r{\rm I}_{r>2}.
\end{align*}
Thus,
\begin{align*}
  &\sum_{k=1}^{\infty} e^{\lambda_0(k-1)}\la R(\theta e^{-\lambda_0k}\phi_0),\psi_0\ra_m\\
  &\le \theta\lambda_0^{-1}/2
    \Big(\int_E\P_{\delta_x}(\theta\la\phi_0,X_1\ra^2;\theta\la\phi_0,X_1\ra\le 2)\psi_0(x)m({\rm d}x)\\
    &\qquad+2C\int_E\P_{\delta_x}(\la\phi_0,X_1\ra
  \log(\theta\la\phi_0,X_1\ra);\theta\la\phi_0,X_1\ra>2)\psi_0(x)m({\rm d}x)\Big).
\end{align*}
Using the dominated convergence theorem, we get
$$\lim_{\theta\to0}\int_E\P_{\delta_x}(\theta\la\phi_0,X_1\ra^2;\theta\la\phi_0,X_1\ra\le 2)\psi_0(x)m({\rm d}x)=0.$$
By Lemma \eqref{lem:llogl}, we have
$$\int_E\P_{\delta_x}(\la\phi_0,X_1\ra\log(\la\phi_0,X_1\ra);\la\phi_0,X_1\ra>2)\psi_0(x)m({\rm d}x)<\infty.$$
Applying the dominated convergence theorem again, we get
$$\lim_{\theta\to0}\int_E\P_{\delta_x}(\la\phi_0,X_1\ra\log(\theta\la\phi_0,X_1\ra);\theta\la\phi_0,X_1\ra>2)\psi_0(x)m({\rm d}x)=0.$$
Therefore, there exists $\theta_0>0$ such that for any
$\theta\in(0, \theta_0]$,
$$\sum_{k=1}^{\infty} e^{\lambda_0(k-1)}\la R(\theta e^{-\lambda_0k}\phi_0),\psi_0\ra_m\le \theta/2.$$
It follows from \eqref{5.15} that, for
$0<\theta\le \theta_0$,
$$\la\Psi(\theta,\cdot),\psi_0\ra_m\ge \theta/2>0.$$
Now the claim \eqref{5.13} is proved,
and hence $M_\infty$ is non-degenerate.

It is easy to see that
\begin{equation}\label{5.19}
 M_\infty=\lim_{t\to\infty}e^{-\lambda_0t}\la\phi_0,X_t\ra=\lim_{t\to\infty}e^{-\lambda_0t}(\gamma_t)^{-1}\gamma_t\la\phi_0,X_t\ra=l_0^{-1}W.
\end{equation}
Since $M_\infty$ is non-degenerate,
 we have $l_0<\infty.$ Thus by \eqref{exp:W},
$\P_{\mu}M_\infty=\la\phi_0,\mu\ra.$
The proof is now complete.
\hfill $\Box$

\noindent
\textbf{Proof of Theorem \ref{llogl3}:}
By \eqref{5.19},
$(1)\Longleftrightarrow(2)\Longleftrightarrow(3)$.
By Proposition \ref{prop1} and Proposition \ref{prop2} , $(3)\Longleftrightarrow(5)\Longleftrightarrow(6)$ and $(4)\Longrightarrow(2)$. Thus, we only need to show that $(3)\Longrightarrow(4)$.

By \eqref{1.8} and the fact that $\eta_s=V_t(\eta_{t+s})$, we have
\begin{align*}
  \eta_0(x)=V_t(\eta_t)(x) = T_t(\eta_t)(x)-\int_0^t T_s\left[r(\cdot,\eta_s(\cdot))\right](x){\rm d}s.
\end{align*}
Thus,
$$
  e^{\lambda_0t}\gamma_t=\la T_t(\eta_t),\psi_0\ra_m=\gamma_0 +\int_E\psi_0(x)m({\rm d}x)\int_0^t e^{\lambda_0s}r(x,\eta_s(x)){\rm d}s,
$$
which implies that
\begin{equation}\label{eq:l}
  l_0=\gamma_0 +\int_E\psi_0(x)m({\rm d}x)\int_0^\infty e^{\lambda_0s}r(x,\eta_s(x)){\rm d}s.
\end{equation}
Recall that,
$$
 r(x,s)=\beta(x)s^2+\int_0^\infty(e^{-s\theta}-1+s\theta)\,n(x,{\rm d}\theta).
$$
Hence,
\begin{align}\label{5.16}
 & \int_0^\infty e^{\lambda_0s}r(x,\eta_s(x)){\rm d}s\ge\int_0^\infty e^{\lambda_0s}{\rm d}s\int_{e^{\lambda_0s}\phi_0(x)^{-1}}^\infty
   (e^{-\theta \eta_s(x)}-1+\theta \eta_s(x))\,n(x,{\rm d}\theta)\nonumber\\
  &=\int_{\phi_0(x)^{-1}}^\infty n(x,{\rm d}\theta )\int_0^{(\log (\theta \phi_0(x))/\lambda_0}e^{\lambda_0s}(e^{-\theta \eta_s(x)}-1+\theta \eta_s(x)){\rm d}s.
\end{align}
Choose $s_0>1$ large enough such that $(1-ce^{-\gamma s_0})>0$,
where $c=c(1)$ and $\gamma=\gamma(1)$ are the constants in \eqref{1.6}.
By \eqref{1.6} and Lemma \ref{lem6}, we get that,
$$
\eta_s(x)=V_{s_0}(\eta_{s+s_0})(x)\ge e^{-as_0}T_{s_0}(\eta_{s+s_0})(x)\ge e^{-as_0}(1-ce^{-\gamma s_0})e^{\lambda_0s_0}\gamma_{s+s_0}\phi_0(x)\ge C\gamma_{s+s_0}\phi_0(x).
$$
Thus,
by Remark \ref{Rem-4.7},
for any $s\le (\log (\theta\phi_0(x))/\lambda_0$, we have
$$
\theta \eta_s(x)\ge C\theta \gamma_{s+s_0}\phi_0(x)\ge
Ce^{\lambda_0s}\gamma_{s+s_0}\ge CL(s+s_0)\ge CL(0)>0,
$$
where $L(t)=e^{\lambda_0 t}\gamma_t$.
Therefore,
$$
\inf_{s\le (\log (\theta \phi_0(x))/\lambda_0}\frac{e^{-\theta \eta_s(x)}-1+\theta \eta_s(x)}{\theta \eta_s(x)}
\ge \inf_{r\ge CL(0)}\frac{e^{-r}-1+r}{r}\ge c,
$$
for some constant $c>0$.
It follows that
\begin{align}\label{I3}
  &\int_{\phi_0(x)^{-1}}^\infty n(x,{\rm d}\theta )\int_0^{(\log (\theta \phi_0(x)))/\lambda_0} e^{\lambda_0s}(e^{-\theta \eta_s(x)}-1+\theta \eta_s(x)){\rm d}s\nonumber \\
   &\ge c\int_{\phi_0(x)^{-1}}^\infty \theta n(x,{\rm d}\theta )\int_0^{\log (\theta\phi_0(x))/\lambda_0} e^{\lambda_0s}\eta_s(x){\rm d}s\nonumber\\
  &\ge C\int_{\phi_0(x)^{-1}}^\infty \theta n(x,{\rm d}\theta )\int_0^{\log (\theta \phi_0(x))/\lambda_0} e^{\lambda_0s}\gamma_{s+s_0}{\rm d}s\phi_0(x)\nonumber\\
  &\ge C\int_{1}^\infty \theta n^{\phi_0}(x,{\rm d}\theta )\int_0^{\log \theta/\lambda_0} L(s+s_0){\rm d}s\nonumber\\
  &\ge CL(0)\int_{1}^\infty \theta\log \theta n^{\phi_0}(x,{\rm d}\theta)\ge Cl(x).
  \end{align}
Combining \eqref{eq:l},\eqref{5.16} and \eqref{I3}, we get that
$l_0\ge C\la l,\psi_0\ra_m.$
Thus $(3)\Longrightarrow(4)$.

The proof is now complete.
\hfill$\Box$

\section{Strong convergence with general test functions}

In this section, we fix a non-trivial family $(\eta_t:t\ge 0)$ of backward iterates and let $\gamma_t
:=\la\eta_t, \psi_0\ra_m$.
The goal of this section is to
determine the almost sure limit of $\gamma_t\langle f, X_t\rangle$
for general test functions $f$.

\subsection{The martingale problems of superprocesses}

In this subsection, we  recall the martingale problem for the superprocess $X$.
Let $J$ denote the set of jump times of $X$, i.e.,
$$
J:=\{s\ge0: \triangle X_s=X_s-X_{s-}\ne0\}.
$$
Since $X$ is a c$\grave{a}$dl$\grave{a}$g process in ${\cal M}_F(E)$, $J$ is a countable set.
 Let $N({\rm d}s,{\rm d}\nu)$ be the optional random measure on
$[0,\infty)\times \mathcal{M}_F(E)$ defined by
$$
N({\rm d}s,{\rm d}\nu):=\sum_{s\in J}\delta_{(s,\triangle X_s)}({\rm d}s,{\rm d}\nu),
$$
and $\widehat{N}({\rm d}s,{\rm d}\nu)$ be the predictable compensator of $N({\rm d}s,{\rm d}\nu)$ which satisfies that
for any nonnegative predictable function $F$ on $\R_+\times
{\cal M}_F(E)\times\Omega$,
\begin{equation}\label{hatN}
\int_0^t\int_{\mathcal{M}_F(E)}F(s,\nu)\,\widehat{N}({\rm d}s,{\rm d}\nu)=\int_0^t{\rm d}s\int_E X_s({\rm d}x)\int_0^\infty
F(s,\theta\delta_x)n(x,{\rm d}\theta),
\end{equation}
where $n$ is the kernel in the branching mechanism $\varphi$.
Define
$$
\widetilde{N}({\rm d}s,{\rm d}\nu):=N({\rm d}s,{\rm d}\nu)-\widehat{N}({\rm d}s,{\rm d}\nu).
$$
Then $\widetilde{N}({\rm d}s,{\rm d}\nu)$ is a martingale measure.
The ``stochastic integral"
$$\int_0^t\int_{\mathcal{M}_F(E)} F(s,\nu)\widetilde{N}({\rm d}s,{\rm d}\nu)
$$
can be defined for a wide class of  Borel functions
$F$  on $[0,t]\times{\cal M}_F(E)$. In particular, if $f$ is a bounded
Borel function on $[0,t]\times E$ and
$F_f(s,\nu):=\int_Ef(s,x)\nu({\rm d}x)$,
then the integral $\int_0^t\int_{\mathcal{M}_F(E)} F_f(s,\nu)\widetilde{N}({\rm d}s,{\rm d}\nu)$ is well-defined.
Let $\mathcal{L}^2_N$ be the space of predictable processes $(F(s,\nu):s>0,\nu\in {\cal M}_F(E))$ satisfying, for all $\mu\in{\cal M}_F(E)$,
$$
\P_{\mu}\int_0^t{\rm d}s\int_E X_s({\rm d}x)\int_0^\infty
F(s,\theta\delta_x)^2n(x,{\rm d}\theta)<\infty.
$$
For any $F\in \mathcal{L}^2_N$,
$$
M^d_t(F):=\int_0^t\int_{\mathcal{M}_F(E)}F(s,\nu)\,\widetilde{N}({\rm d}s,{\rm d}\nu),\qquad t\ge0,
$$
is a square integrable martingale such that
\begin{equation}\label{sqr-N}
  \P_{\mu}(M^d_t(F)^2)=\P_{\mu}\left[\int_0^t{\rm d}s\int_E X_s({\rm d}x)\int_0^\infty
F(s,\theta\delta_x)^2n(x,{\rm d}\theta)\right].
\end{equation}

Note that $C_0(E)$ is a Banach space under the supremum norm.
In the remainder of this paper, we assume that
\begin{assumption}\label{assum5a}
\begin{description}
 \item[(i)]
 $\{P_t, t\ge 0\}$ is a Feller semigroup, that is, $\{P_t, t\ge 0\}$ preserves $C_0(E)$ and
 $[0, \infty) \ni t\rightarrow P_tf\in C_0(E)$ is continuous  for every $f\in C_0(E).$

 \item[(ii)]
\begin{equation}\label{limit-n-0}
\lim_{a\to\infty}\sup_{x\in E}\int_a^\infty \theta n(x,{\rm d}\theta)=0.
\end{equation}
\end{description}
\end{assumption}

In the reminder of this subsection, we will (also) use $\widetilde{L}$ to denote the infinitesimal generator of $\{P_t, t\ge 0\}$ in the space $C_0(E)$ and use ${\rm Dom}(\widetilde{L})$ its domain.
It is known (see \cite[Section 6.1]{Dawson}, for instance) that,
for any $f\in {\rm Dom}(\widetilde{L})$, we have that
$$
  \la f,X_t\ra=\la f,X_0\ra+\int_0^t \la
\widetilde{L}f-\alpha f,
X_s\ra{\rm d}s+M_t^c(f)+M_t^d(f),
$$
where $M_t^c(f)$ is a continuous local martingale with quadratic variation
$\int^t_02\la \beta f^2,X_s\ra{\rm d}s$ and
$$
M_t^d(f)=\int_0^t\int_{\mathcal{M}_F(E)} \la f,\nu\ra\widetilde{N}({\rm d}s,{\rm d}\nu)
$$
is a purely discontinuous local martingale.
Here we remark that if we assume  $\alpha, \beta\in C(E)$ and that
$x\to (\theta\wedge \theta^2)n(x, {\rm d}\theta)$
is continuous in the topology of weak convergence, then the above result follows from \cite[Theorem 7.25]{Li11}.
$M_t^c(f)$  induces a  worthy $({\cal G}_t)$-martingale measure $S^C({\rm d}s,{\rm d}x)$
(see \cite[Section 7.3]{Li11} for the definition of worthy martingale measure) satisfying
$$
M_t^c(f)=\int_0^t\int_E f(x)S^C({\rm d}s,{\rm d}x).
$$
Standard arguments then show that the ``stochastic integral"
$$
\int_0^t\int_E f(s,x)S^C({\rm d}s,{\rm d}x)
$$
can be defined for a wide class of integrands $f$ on $[0,t]\times E$ (see, for example,  \cite[Theorem 7.25]{Li11} or \cite{Fit} for more details).
Thus, one can show that (see \cite[Proposition 2.13]{Fit} or \cite[Exercise II.5.2]{Perkins} for
the case when the branching mechanism has finite variance)
for any bounded function $g$ on $E$,
\begin{equation}\label{martep2}
  \la g,X_t\ra=\la T_tg, X_0 \ra+\int_0^t\int_{\mathcal{M}_F(E)}\la T_{t-s}g,\nu\ra \widetilde{N}({\rm d}s,{\rm d}\nu)+\int_0^t\int_ET_{t-s}g(x)S^C({\rm d}s,{\rm d}x).
\end{equation}

\subsection{Discrete times}
In this subsection,
we show that for any $\delta>0$ and $f\in \mathcal{B}_b^+(E)$,  $\gamma_{n\delta}\langle f\phi_0, X_{n\delta}\rangle$ has an almost sure limit as $n\to\infty$.
 We will extend this result to continuous times in two different scenarios in the next two subsections.

\begin{theorem}\label{them:dis}
For any $\delta>0$, $\mu\in \mathcal{M}_F(E)$ and $f\in \mathcal{B}_b^+(E)$, we have
$$
\lim_{n\to\infty}\gamma_{n\delta}\la \phi_0 f,X_{n\delta}\ra
=\la f\phi_0,\psi_0\ra_m W,\qquad \mbox{a.s.-}\P_{\mu}.
$$
\end{theorem}

To prove Theorem \ref{them:dis}, we first make some preparations.
For any $s>0,$  we define
\begin{align}\label{def:D1}
\mathcal{D}_{<1}(s):
=\{\nu\in \mathcal{M}_F(E):0<\gamma_s\la\phi_0,\nu\ra<1\}
\end{align}
and
\begin{align}\label{def:D2}
\mathcal{D}_{\ge1}(s):
=\{\nu\in \mathcal{M}_F(E):\gamma_s\la\phi_0,\nu\ra\ge1\}.
\end{align}
For any $m\in \mathbb{N}$, $\delta>0$, $\mu\in \mathcal{M}_F(E)$ and $f\in \mathcal{B}_b^+(E)$, by \eqref{martep2}, we have
\begin{align*}
&  \gamma_{(n+m)\delta}\la\phi_0f,X_{(n+m)\delta}\ra\\
  &=\gamma_{(n+m)\delta}\la T_{(n+m)\delta}(\phi_0f), \mu \ra
  +\gamma_{(n+m)\delta}\int_0^{(n+m)\delta}\int_ET_{(n+m)\delta-s}(\phi_0f)(x)S^C({\rm d}s,{\rm d}x)\\
  &\qquad+\gamma_{(n+m)\delta}\int_0^{(n+m)\delta}
  \int_{\mathcal{D}_{<1}(s)}
  \la T_{(n+m)\delta-s}(\phi_0f),\nu\ra \widetilde{N}({\rm d}s,{\rm d}\nu) \\
   &\qquad+\gamma_{(n+m)\delta}\int_0^{(n+m)\delta}
   \int_{\mathcal{D}_{\ge 1}(s)}
   \la T_{(n+m)\delta-s}(\phi_0f),\nu\ra \widetilde{N}({\rm d}s,{\rm d}\nu) \\
   &=: \gamma_{(n+m)\delta}\la T_{(n+m)\delta}(\phi_0f), \mu \ra
   +C_{(n+m)\delta}(f)+H_{(n+m)\delta}(f)+L_{(n+m)\delta}(f) .
\end{align*}
Therefore,
\begin{align*}
  &\gamma_{(n+m)\delta}\la\phi_0f,X_{(n+m)\delta}\ra
  -\P_{\mu}[\gamma_{(n+m)\delta}\la\phi_0f,X_{(n+m)\delta}\ra|\mathcal{G}_{n\delta}] \\
   &=  \left(H_{(n+m)\delta}(f)-\P_{\mu}(H_{(n+m)\delta}(f)|\mathcal{G}_{n\delta})\right)
   +\left(L_{(n+m)\delta}(f)-\P_{\mu}(L_{(n+m)\delta}(f)|\mathcal{G}_{n\delta})\right)\\
   &\qquad+\left(C_{(n+m)\delta}(f)-\P_{\mu}(C_{(n+m)\delta}(f)|\mathcal{G}_{n\delta})\right).
\end{align*}
We now deal with the three  parts separately.
Before doing this, we prove a lemma first.

\begin{lemma}\label{lem:3.1}
If $\{a_n: n\ge 1\}$ is a sequence of positive numbers such that
$\lim_{n\to\infty}a_{n+1}/a_n=a>1$, then
\begin{equation}\label{4.7}
  \sup_{x\in E}\sum_{n=1}^\infty a_n^{-1}\int_0^{a_n}
\theta^2 n^{\phi_0}(x,{\rm d}\theta)<\infty
\end{equation}
and
\begin{equation}\label{4.8}
  \sup_{x\in E}\sum_{n=1}^\infty a_n\int_{a_n}^\infty
n^{\phi_0}(x,{\rm d}\theta)<\infty.
\end{equation}
\end{lemma}

\noindent\textbf{Proof:}
Since
$\lim_{n\to\infty} a_{n+1}/a_n=a>1$,  for any $a^*\in(1,a)$,
there exists $K>0$ such that for any $n\ge K$,
$$
\frac{a_{n+1}}{a_n}>a^*.
$$
Without loss of generality, we assume that $a_n\uparrow\infty$.
For convenience, we put $a_0=0$.
Then we have
\begin{align*}
  &\sum_{n=1}^\infty a_n^{-1}\int_0^{a_n} \theta^2 n^{\phi_0}(x,{\rm d}\theta)
  =\sum_{n=1}^\infty a_n^{-1}\sum_{k=1}^{n}\int_{a_{k-1}}^{a_k} \theta^2 n^{\phi_0}(x,{\rm d}\theta)\\
   &= \sum_{k=1}^\infty \int_{a_{k-1}}^{a_k} \theta^2 n^{\phi_0}(x,{\rm d}\theta)\sum_{n=k}^{\infty}a_n^{-1}\\
   &\le \sum_{k=1}^K \int_{a_{k-1}}^{a_k} \theta^2 n^{\phi_0}(x,{\rm d}\theta) \sum_{n=k}^{\infty}a_n^{-1}
    +\sum_{k=K+1}^\infty \int_{a_{k-1}}^{a_k} \theta n^{\phi_0}(x,{\rm d}\theta) a_k\sum_{n=k}^{\infty}a_n^{-1}\\
   &\le  \int_{0}^{a_K} \theta^2 n^{\phi_0}(x,{\rm d}\theta)\sum_{n=1}^{\infty}a_n^{-1}
   +\sum_{k=K+1}^\infty \int_{a_{k-1}}^{a_k} \theta n^{\phi_0}(x,{\rm d}\theta) a_k\sum_{n=k}^{\infty}a_n^{-1}.
\end{align*}
For any $k>K$, we have
$$
\sum_{n=k}^{\infty}a_ka_n^{-1}\le \sum_{n=k}^{\infty}(a^*)^{-(n-k)}=\sum_{n=1}^{\infty}(a^*)^{-n}<\infty.
$$
Therefore,
\begin{align*}
  &\sum_{n=1}^\infty a_n^{-1}\int_0^{a_n} \theta^2 n^{\phi_0}(x,{\rm d}\theta)
   \le C\left[\int_{0}^{a_K} \theta^2 n^{\phi_0}(x,{\rm d}\theta)+\int_{a_K}^{\infty} \theta n^{\phi_0}(x,{\rm d}\theta)\right]\\
   &\le  C\sup_{x\in E}\int_{0}^{\infty}
   (\theta\wedge \theta^2) n^{\phi_0}(x,{\rm d}\theta)<\infty.
\end{align*}

Note that
\begin{align*}
  &\sum_{n=1}^\infty a_n\int_{a_n}^\infty  n^{\phi_0}(x,{\rm d}\theta)
  =\sum_{n=1}^\infty a_n\sum_{k=n}^{\infty}
  \int_{a_{k}}^{a_{k+1}}  n^{\phi_0}(x,{\rm d}\theta)\\
   &\le  \sum_{k=1}^\infty
    \int_{a_{k}}^{a_{k+1}} \theta n^{\phi_0}(x,{\rm d}\theta)\Big(a_k^{-1}\sum_{n=1}^{k}a_n\Big).
\end{align*}
Using elementary calculus, one can easily show that
$$
\lim_{k\to\infty}a_k^{-1}\sum_{n=1}^{k}a_n=\lim_{k\to\infty}\frac{a_{k+1}}{a_{k+1}-a_k}
=\frac{a}{a-1}.
$$
Thus $\sup_{k\ge 1}a_k^{-1}\sum_{n=1}^{k}a_n<\infty$.
It follows that
$$
\sup_{x\in E}\sum_{n=1}^\infty a_n\int_{a_n}^\infty
n^{\phi_0}(x,{\rm d}\theta)\le C\sup_{x\in E}\int_{a_1}^{\infty} \theta n^{\phi_0}(x,{\rm d}\theta)
<\infty.
$$
The proof is now complete.

\hfill $\Box$

Define
$$
I(a,x):=\int_0^a \theta^2n^{\phi_0}(x,{\rm d}\theta).
$$

\begin{lemma}\label{lem:H}
For any $m\in \mathbb{N}$, $\delta>0$, $\mu\in \mathcal{M}_F(E)$ and $f\in \mathcal{B}_b^+(E)$, we have
$$\lim_{n\to\infty}H_{(n+m)\delta}(f)-\P_{\mu}(H_{(n+m)\delta}(f)|\mathcal{G}_{n\delta})=0,
\qquad \mbox{a.s.-}\P_{\mu}.$$
\end{lemma}

\noindent\textbf{Proof:}
By the conditional Borel-Cantelli lemma, it suffices to prove that
\begin{equation}\label{sqr:H}
  \sum_{n=1}^\infty\P_{\mu}\left([H_{(n+m)\delta}(f)
  -\P_{\mu}(H_{(n+m)\delta}(f)|\mathcal{G}_{n\delta})]^2|\mathcal{G}_{(n-1)\delta}\right)<\infty.
\end{equation}
Recall from (\ref{def:D1}) that
\begin{align*}
\mathcal{D}_{<1}(s):=
\{\nu\in \mathcal{M}_F(E):0<\gamma_s\la\phi_0,\nu\ra<1\}.
\end{align*}
Since $\widetilde{N}({\rm d}s,{\rm d}\nu)$ is a martingale measure, we have
$$\P_{\mu}(H_{(n+m)\delta}(f)|\mathcal{G}_{n\delta})
=\gamma_{(n+m)\delta}\int_0^{n\delta}
\int_{\mathcal{D}_{<1}(s)}
\la T_{(n+m)\delta-s}(\phi_0f),\nu\ra \widetilde{N}({\rm d}s,{\rm d}\nu),$$
which implies that
\begin{equation}\label{est-H}
  H_{(n+m)\delta}(f)-\P_{\mu}(H_{(n+m)\delta}(f)|\mathcal{G}_{n\delta})
  =\gamma_{(n+m)\delta}\int_{n\delta}^{(n+m)\delta}
  \int_{\mathcal{D}_{<1}(s)}
  \la T_{(n+m)\delta-s}(\phi_0f),\nu\ra \widetilde{N}({\rm d}s,{\rm d}\nu).
\end{equation}
By \eqref{eq:eta},
$\gamma_t=\la\Phi(e^{-\lambda_0t},\cdot),\psi_0\ra_m$,
which implies that $t\to \gamma_t$ is non-increasing.
Thus by \eqref{sqr-N} and \eqref{1.6}, we have
\begin{align}\label{4.33}
  &\P_{\mu}\left([H_{(n+m)\delta}(f)-\P_{\mu}(H_{(n+m)\delta}(f)|\mathcal{G}_{n\delta})]^2|\mathcal{G}_{(n-1)\delta}\right)\nonumber\\
  &=\gamma_{(n+m)\delta}^2\P_{X_{(n-1)\delta}}\left[\int_{\delta}^{(1+m)\delta}{\rm d}s
  \int_E X_s({\rm d}x)\int_0^{\gamma_{s+(n-1)\delta}^{-1}\phi_0(x)^{-1}}
\theta^2[T_{(1+m)\delta-s}(\phi_0f)(x)]^2n(x,{\rm d}\theta)\right]\nonumber\\
  &\le\|f\|_\infty^2 \gamma^2_{(n+m)\delta}\P_{X_{(n-1)\delta}}\left[\int_{\delta}^{(1+m)\delta} e^{2\lambda_0((1+m)\delta-s)}{\rm d}s
  \int_E X_s({\rm d}x)\int_0^{\gamma_{(n+m)\delta}^{-1}}
\theta^2n^{\phi_0}(x,{\rm d}\theta)\right]\nonumber\\
  &=\|f\|_\infty^2 \gamma^2_{(n+m)\delta}\int_{\delta}^{(1+m)\delta} e^{2\lambda_0((1+m)\delta-s)}
  \la T_s[I(\gamma_{(n+m)\delta}^{-1})],X_{(n-1)\delta}\ra{\rm d}s\nonumber\\
  &\le(1+c(\delta))\|f\|_\infty^2 e^{2\lambda_0(1+m)\delta}\int_{\delta}^{(1+m)\delta} e^{-\lambda_0s}{\rm d}s\gamma^2_{(n+m)\delta}
  \la I(\gamma_{(n+m)\delta}^{-1}),\psi_0\ra_m\la \phi_0,X_{(n-1)\delta}\ra\nonumber\\
  &\le C\left[\gamma_{(n+m)\delta}\la \phi_0,X_{(n-1)\delta}\ra \right]
  \gamma_{(n+m)\delta} \la I(\gamma_{(m+n)\delta}^{-1}),\psi_0\ra_m.
\end{align}
It follows from the fact that
$\lim_{n\to\infty}\frac{\gamma_{(n+m)\delta}^{-1}}{\gamma_{(n+m-1)\delta}^{-1}}=e^{\lambda_0\delta}$
and Lemma \ref{lem:3.1} that
\begin{equation}\label{4.1}
  \sum_{n=1}^\infty\gamma_{(n+m)\delta} \la I(\gamma_{(m+n)\delta}^{-1}),\psi_0\ra_m<\infty.
\end{equation}
Since $\lim_{n\to\infty}\gamma_{(n+m)\delta}\la \phi_0,X_{(n-1)\delta}\ra=e^{-\lambda_0(m+1)\delta}W$,
combining \eqref{4.33} and \eqref{4.1},
\eqref{sqr:H} follows immediately.
The proof is now complete.
\hfill $\Box$

\begin{lemma}\label{lem:L}
For any $m\in\mathbb{N}$, $\delta>0$, $\mu\in \mathcal{M}_F(E)$ and $f\in \mathcal{B}_b^+(E)$, we have
$$
\lim_{n\to\infty}L_{(n+m)\delta}(f)-\P_{\mu}(L_{(n+m)\delta}(f)|\mathcal{G}_{n\delta})=0,
\qquad \mbox{a.s.-}\P_{\mu}.
$$
\end{lemma}
\noindent\textbf{Proof:}
Recall the definition of
$\mathcal{D}_{\ge 1}(s)$ in (\ref{def:D2}):
$$
\mathcal{D}_{\ge 1}(s)=
\{\nu\in \mathcal{M}_F(E):\gamma_s\la\phi_0,\nu\ra\ge1\}.
$$
Since $\widetilde{N}({\rm d}s,{\rm d}\nu)$ is a martingale measure, we have
$$\P_{\mu}(L_{(n+m)\delta}(f)|\mathcal{G}_{n\delta})
=\gamma_{(n+m)\delta}\int_0^{n\delta}
\int_{\mathcal{D}_{\ge 1}(s)}
\la T_{(n+m)\delta-s}(\phi_0f),\nu\ra \widetilde{N}({\rm d}s,{\rm d}\nu),$$
which implies that
\begin{align*}
  &L_{(n+m)\delta}(f)-\P_{\mu}(L_{(n+m)\delta}(f)|\mathcal{G}_{n\delta})\\
  &=\gamma_{(n+m)\delta}\int_{n\delta}^{(n+m)\delta}
  \int_{\mathcal{D}_{\ge 1}(s)}
  \la T_{(n+m)\delta-s}(\phi_0f),\nu\ra \Big(N({\rm d}s,{\rm d}\nu)-\widehat{N}({\rm d}s,{\rm d}\nu\Big).
\end{align*}

We claim that
\begin{equation}\label{4.4}
  \P_{\mu}\Big(\int_{n\delta}^{(n+m)\delta}
  \int_{\mathcal{D}_{\ge 1}(s)}N({\rm d}s,{\rm d}\nu)>0,\,\mbox{i.o.}\Big)=0.
\end{equation}
In fact, since
$\int_{n\delta}^{(n+m)\delta}
\int_{\mathcal{D}_{\ge 1}(s)}N({\rm d}s,{\rm d}\nu)$ is a non-negative integer,
by the Markov property of $X$,
\begin{align}\label{4.6}
  &\sum_{n=1}^\infty\P_{\mu}\Big(\int_{n\delta}^{(n+m)\delta}
  \int_{\mathcal{D}_{\ge 1}(s)}
  N({\rm d}s,{\rm d}\nu)>0\Big|\mathcal{G}_{(n-1)\delta}\Big)\nonumber\\
 &\le \sum_{n=1}^\infty\P_{\mu}\Big(\int_{n\delta}^{(n+m)\delta}
 \int_{\mathcal{D}_{\ge 1}(s)}
 N({\rm d}s,{\rm d}\nu)\Big|\mathcal{G}_{(n-1)\delta}\Big)\nonumber\\
 &=\sum_{n=1}^\infty\P_{X_{(n-1)\delta}}\Big(\int_{\delta}^{(1+m)\delta}
 \int_{\mathcal{D}_{\ge 1}(s+(n-1)\delta)}
 N({\rm d}s,{\rm d}\nu)\Big)\nonumber\\
 &=\sum_{n=1}^\infty\P_{X_{(n-1)\delta}}
 \left[\int_{\delta}^{(1+m)\delta}ds\int_{E}X_s({\rm d}x)
 \int_{\phi_0(x)^{-1}\gamma_{s+(n-1)\delta}^{-1}}^\infty
n(x,{\rm d}\theta)\right]\nonumber\\
 &\le\sum_{n=1}^\infty\int_{\delta}^{(1+m)\delta}{\rm d}s \left
 \la T_s\Big[\int_{\gamma_{n\delta}^{-1}}^\infty
n^{\phi_0}(\cdot,{\rm d}\theta)\Big],
 X_{(n-1)\delta}\right\ra\nonumber\\
 &\le(1+c(\delta))m\delta e^{\lambda_0(m+1)\delta}\sum_{n=1}^\infty
 \left\la \int_{\gamma_{n\delta}^{-1}}^\infty
n^{\phi_0}(\cdot,{\rm d}\theta),\psi_0\right\ra_m
 \la \phi_0,X_{(n-1)\delta}\ra,
\end{align}
where in the  second to the last inequality, we use the fact that
$\gamma_{s+(n-1)\delta}\le \gamma_{n\delta}$,
and the last inequality follows from \eqref{1.6}.
It follows from \eqref{4.8} that
\begin{equation}\label{4.9}
\sum_{n=1}^\infty \gamma_{n\delta}^{-1}\la \int_{\gamma_{n\delta}^{-1}}^\infty
n^{\phi_0}(\cdot,{\rm d}\theta),\psi_0\ra_m<\infty.
\end{equation}
By Theorem \ref{Thm:main1},
$\gamma_{n\delta}\la \phi_0,X_{(n-1)\delta}\ra\to e^{-\lambda_0\delta}W$ as $n\to\infty$.
Therefore we have
$$
\sum_{n=1}^\infty \la \int_{\gamma_{n\delta}^{-1}}^\infty
n^{\phi_0}(\cdot,{\rm d}\theta),\psi_0\ra_m
\la \phi_0,X_{(n-1)\delta}\ra<\infty,
$$
which implies that
$$
\sum_{n=1}^\infty\P_{\mu}\Big(\int_{n\delta}^{(n+m)\delta}
\int_{\mathcal{D}_{\ge 1}(s)}
N({\rm d}s,{\rm d}\nu)>0\Big|\mathcal{G}_{(n-1)\delta}\Big)<\infty.
$$
Now using the the conditional Borel-Cantelli lemma, we immediately get
the claim \eqref{4.4}.

By \eqref{4.4}, we get
\begin{equation}\label{4.5}
  \lim_{n\to\infty}\gamma_{(n+m)\delta}\int_{n\delta}^{(n+m)\delta}
 \int_{\mathcal{D}_{\ge 1}(s)}
  \la T_{(n+m)\delta-s}(\phi_0f),\nu\ra N({\rm d}s,{\rm d}\nu)=0,\qquad \P_{\mu}\mbox{-a.s.}
\end{equation}
To complete the proof, we only need to show that
\begin{equation}\label{4.10}
  \lim_{n\to\infty}\gamma_{(n+m)\delta}\int_{n\delta}^{(n+m)\delta}
  \int_{\mathcal{D}_{\ge 1}(s)}
  \la T_{(n+m)\delta-s}(\phi_0f),\nu\ra \widehat{N}({\rm d}s,{\rm d}\nu)=0,\qquad \P_{\mu}\mbox{-a.s.}
\end{equation}
By \eqref{hatN}, we have
\begin{align*}
  &\gamma_{(n+m)\delta}\int_{n\delta}^{(n+m)\delta}
  \int_{\mathcal{D}_{\ge 1}(s)}
  \la T_{(n+m)\delta-s}(\phi_0f),\nu\ra \widehat{N}({\rm d}s,{\rm d}\nu) \\
   &= \gamma_{(n+m)\delta}\int_{n\delta}^{(n+m)\delta}{\rm d}s\int_{E}T_{(n+m)\delta-s}(\phi_0f)(x)X_s({\rm d}x)
   \int_{\phi_0(x)^{-1}\gamma_s^{-1}}^\infty
\theta n(x,{\rm d}\theta)\\
   &\le\|f\|_\infty e^{\lambda_0m\delta}\gamma_{(n+m)\delta}
   \int_{n\delta}^{(n+m)\delta}{\rm d}s\int_{E}\phi_0(x)X_s({\rm d}x)
   \int_{\phi_0(x)^{-1}\gamma_s^{-1}}^\infty
\theta n(x,{\rm d}\theta)\\
   &\le\|f\|_\infty e^{\lambda_0m\delta}\gamma_{(n+m)\delta}
   \int_{n\delta}^{(n+m)\delta}\la\phi_0,X_s\ra {\rm d}s\sup_{x\in E}
   \left(\int_{\|\phi_0\|_\infty^{-1}\gamma_{n\delta}^{-1}}^\infty
\theta n(x,{\rm d}\theta)\right).
\end{align*}
Note that
\begin{align*}
  &\lim_{n\to\infty}\gamma_{(n+m)\delta}\int_{n\delta}^{(n+m)\delta}\la\phi_0,X_s\ra {\rm d}s =\lim_{n\to\infty}\int_{0}^{m\delta}\gamma_{(n+m)\delta}\la\phi_0,X_{s+n\delta}\ra {\rm d}s\\
  &=\int_{0}^{m\delta} e^{-\lambda_0(m\delta-s)}{\rm d}s W,
\end{align*}
 and by \eqref{limit-n-0} we have
 $$
 \lim_{n\to\infty}\sup_{x\in E}
 \left(\int_{\|\phi_0\|_\infty^{-1}\gamma_{n\delta}^{-1}}^\infty
\theta n(x,{\rm d}\theta)\right)=0.
 $$
Now we easily see that \eqref{4.10} holds. The proof is now complete.
\hfill $\Box$

\begin{lemma}\label{lem:C}
For any $m\in \mathbb{N}$, $\delta>0$, $\mu\in \mathcal{M}_F(E)$ and $f\in \mathcal{B}_b^+(E)$, we have
\begin{equation}\label{lim:C}
  \lim_{n\to\infty}C_{(n+m)\delta}(f)-\P_{\mu}(C_{(n+m)\delta}(f)|\mathcal{G}_{n\delta})=0,
  \qquad \mbox{a.s.-}\P_{\mu}.
\end{equation}
\end{lemma}
\noindent\textbf{Proof:}
Let
$$
\widetilde{M}_t:=\int_0^{t}\int_E T_{(n+m)\delta-s}(\phi_0f)(x)S^C({\rm d}s,{\rm d}x).
$$
Then $\{\widetilde{M}_t,0\le t\le (n+m)\delta\}$ is a martingale with quadratic variation
$$
\la \widetilde{M}\ra_t=2\int_0^{t}\la\beta( T_{(n+m)\delta-s}(\phi_0f))^2,X_s\ra{\rm d}s.
$$
Note that
\begin{align}\label{est-C}
  &C_{(n+m)\delta}(f)-\P_{\mu}(C_{(n+m)\delta}(f)|\mathcal{G}_{n\delta})\nonumber\\
  &=\gamma_{(n+m)\delta}\int_{n\delta}^{(n+m)\delta}
  \int_{E} T_{(n+m)\delta-s}(\phi_0f)(x) S^C({\rm d}s,{\rm d}x)\nonumber\\
  &=\gamma_{(m+n)\delta}\,(\widetilde{M}_{(n+m)\delta}-\widetilde{M}_{n\delta}).
\end{align}
Using this we get
\begin{align*}
  &\P_{\mu}\left([C_{(n+m)\delta}(f)-\P_{\mu}(C_{(n+m)\delta}(f)|\mathcal{G}_{n\delta})]^2|\mathcal{G}_{(n-1)\delta}\right)\nonumber\\
  &=\gamma_{(m+n)\delta}^2\P_{\mu}\left(\la \widetilde{M}\ra_{(n+m)\delta}-\la \widetilde{M}\ra_{n\delta}|\mathcal{G}_{(n-1)\delta}\right)\\
  &=\gamma_{(m+n)\delta}^2\P_{\mu}\left(2\int_{n\delta}^{(n+m)\delta}\la \beta( T_{(n+m)\delta-s}(\phi_0f))^2,X_s\ra
  {\rm d}s|\mathcal{G}_{(n-1)\delta}\right)\\
  &=\gamma_{(m+n)\delta}^2\P_{X_{(n-1)\delta}}\left(2\int_{\delta}^{(1+m)\delta}\la \beta( T_{(1+m)\delta-s}(\phi_0f))^2,X_s\ra{\rm d}s\right).
\end{align*}
Note that
$$
\beta(x)( T_{(1+m)\delta-s}(\phi_0f))^2(x)
\le \|f\|_\infty^2\beta(x)e^{2\lambda_0( (1+m)\delta-s)}\phi_0(x)^2
\le \|f\|_\infty^2\|\beta  \phi_0\|_\infty e^{2\lambda_0 m\delta}\phi_0(x).
$$
Thus we have
\begin{align*}
   &\P_{\mu}\left([C_{(n+m)\delta}(f)-\P_{\mu}(C_{(n+m)\delta}(f)|\mathcal{G}_{n\delta})]^2|\mathcal{G}_{(n-1)\delta}\right)\nonumber\\
   &\le C\gamma_{(m+n)\delta}^2\P_{X_{(n-1)\delta}}\int_{\delta}^{(1+m)\delta}\la \phi_0,X_s\ra
    {\rm d}s=C\gamma^2_{(m+n)\delta}\int_{\delta}^{(1+m)\delta}\la T_s\phi_0,X_{(n-1)\delta}\ra{\rm d}s\\
   &=C\gamma_{(m+n)\delta}^2\int_{\delta}^{(1+m)\delta}e^{\lambda_0 s}{\rm d}s\,\la \phi_0,X_{(n-1)\delta}\ra
   \le C\gamma_{(m+n)\delta}\Big(\gamma_{(m+n)\delta}\la \phi_0,X_{(n-1)\delta}\ra\Big).
 \end{align*}
By \eqref{3.10}, we have that
$\lim_{n\to\infty}\frac{\gamma_{(m+n)\delta}}{\gamma_{(m+n-1)\delta}}=e^{-\lambda_0\delta}<1$,
which implies that
$$
\sum_{n=1}^\infty \gamma_{(m+n)\delta}<\infty.
$$
By Theorem \ref{Thm:main1},
$$
\lim_{n\to\infty}\gamma_{(m+n)\delta}\la \phi_0,X_{(n-1)\delta}\ra=e^{-\lambda_0(m+1)\delta}W.
$$
Thus we have
$$\sum_{n=1}^\infty\gamma_{(m+n)\delta}\Big(\gamma_{(m+n)\delta}\la \phi_0,X_{(n-1)\delta}\ra\Big)<\infty,$$
which implies that
$$\sum_{n=1}^\infty\P_{\mu}\left([C_{(n+m)\delta}(f)-\P_{\mu}(C_{(n+m)\delta}(f)|\mathcal{G}_{n\delta})]^2|
\mathcal{G}_{(n-1)\delta}\right)<\infty.$$
Now  using the conditional Bore-Cantelli lemma, we immediately get \eqref{lim:C}.
\hfill $\Box$

Combining the three results above, we get

\begin{lemma}\label{lem:cod}
For any $m\in \mathbb{N}$, $\delta>0$, $\mu\in \mathcal{M}_F(E)$ and $f\in \mathcal{B}_b^+(E)$,
we have
$$\lim_{n\to\infty} \gamma_{(n+m)\delta}\la\phi_0f,X_{(n+m)\delta}\ra
-\P_{\mu}[\gamma_{(n+m)\delta}\la\phi_0f,X_{(n+m)\delta}\ra|\mathcal{G}_{n\delta}]=0,\qquad \mbox{a.s.-}\P_{\mu}.
$$
\end{lemma}

\bigskip
\noindent
\textbf{Proof of Theorem \ref{them:dis}:}
By the Markov property of $X$, we have
$$
\gamma_{(n+m)\delta}\P_{\mu}(\la \phi_0 f,X_{(n+m)\delta}\ra|\mathcal{G}_{n\delta})
=\gamma_{(n+m)\delta}\la T_{m\delta}(\phi_0 f),X_{n\delta}\ra.
$$
It follows from \eqref{T_t} that there exist constants $c>0$ and $\gamma>0$ such that for any $m\ge1$,
$$
(1-ce^{-\gamma m\delta})e^{\lambda_0m\delta}\la f\phi_0,\psi_0\ra_m\phi_0(x)\le T_{m\delta}(\phi_0 f)(x)\le (1+ce^{-\gamma m\delta})e^{\lambda_0m\delta}\la f\phi_0,\psi_0\ra_m\phi_0(x).
$$
Thus, by Lemma \ref{lem:cod}, we have that, for any $m\in \mathbb{N}$,
\begin{align*}
  \limsup_{n\to\infty} \gamma_{n\delta}\la \phi_0 f,X_{n\delta}\ra
  &= \limsup_{n\to\infty} \gamma_{(n+m)\delta}\la \phi_0 f,X_{(n+m)\delta}\ra \\
   &= \limsup_{n\to\infty} \gamma_{(n+m)\delta}\la T_{m\delta}(\phi_0 f),X_{n\delta}\ra\\
   &\le \limsup_{n\to\infty} \gamma_{(n+m)\delta}(1+ce^{-\gamma m\delta})
   e^{\lambda_0m\delta}\la f\phi_0,\psi_0\ra_m\la \phi_0,X_{n\delta}\ra\\
   &=(1+ce^{-\gamma m\delta})\la f\phi_0,\psi_0\ra_m W.
\end{align*}
Letting $m\to\infty$, we get
\begin{equation}\label{limsup}
   \limsup_{n\to\infty} \gamma_{n\delta}\la \phi_0 f,X_{n\delta}\ra\le \la f\phi_0,\psi_0\ra_m W.
\end{equation}
Similarly, we have
$$ \liminf_{n\to\infty} \gamma_{n\delta}\la \phi_0 f,X_{n\delta}\ra\ge(1-ce^{-\gamma m\delta})\la f\phi_0,\psi_0\ra_m W.$$
Letting $m\to\infty$, we get
\begin{equation}\label{liminf}
   \liminf_{n\to\infty} \gamma_{n\delta}\la \phi_0 f,X_{n\delta}\ra\ge \la f\phi_0,\psi_0\ra_m W.
\end{equation}
 Combining \eqref{limsup} and \eqref{liminf}, the conclusion follows immediately.
 \hfill$\Box$

\subsection{Continuous times: Case I}
Define a new semigroup $(T^{\phi_0}_t,t\ge0)$ by
 $$
 T_t^{\phi_0}f(x):=\frac{e^{-\lambda_0t}T_t(f\phi_0)(x)}{\phi_0(x)},\qquad
 f\in \mathcal{B}_b(E).
 $$
Then $(T_t^{\phi_0},t\ge0)$ is conservative semigroup
with transition density
 $$
 q^{\phi_0}(t,x,y)=\frac{e^{-\lambda_0t}q(t,x,y)\phi_0(y)}{\phi_0(x)}.
 $$

In this subsection, we also make the following assumption:

\begin{assumption}\label{assum5}
For any $f\in C_0(E)$,
\begin{equation}\label{Feller}\lim_{t\to 0}\|T^{\phi_0}_tf-f\|_\infty=0.\end{equation}
\end{assumption}

See Examples 4.4, 4.5,  4.7 and Remark 4.6 of \cite{CRSZ} for examples satisfying the assumption above,
and Assumptions 1 and 2.

\begin{theorem}\label{thrm:conti}
Under Assumptions  \ref{assum1}--\ref{assum5},
we have that, for any $\mu\in \mathcal{M}_F(E)$ and
$f\in C_0(E)$,
\begin{equation}\label{conti1}
\lim_{t\to\infty}\gamma_t\la \phi_0f,X_t\ra=\la f\phi_0,\psi_0\ra_m W,\qquad \mbox{a.s.-}\P_{\mu}.
\end{equation}
\end{theorem}
\noindent\textbf{Proof:}
First, we claim that
\begin{equation}\label{claim1}
\lim_{\delta\to0}\limsup_{n\to\infty} \sup_{t\in[n\delta,(n+1)\delta]}
\left|\gamma_t\la \phi_0T^{\phi_0}_{(n+1)\delta-t}f,X_t\ra- \gamma_t\la \phi_0f,X_t\ra\right|=0.
\end{equation}
In fact, we have
\begin{align*}
   &\sup_{t\in[n\delta,(n+1)\delta]}\left|\gamma_t\la \phi_0T^{\phi_0}_{(n+1)\delta-t}f,X_t\ra- \gamma_t\la \phi_0f,X_t\ra\right|\\
  &\le   \sup_{r\in(0,\delta)}\|T^{\phi_0}_rf-f\|_\infty\sup_{t\in[n\delta,(n+1)\delta]}\gamma_t\la\phi_0,X_t\ra.
\end{align*}
Letting $n\to\infty$ and then $\delta\to0$, using Assumption \ref{assum5} and Theorem \ref{Thm:main1}, one immediately arrives at
the claim \eqref{claim1}.
Thus, by \eqref{claim1}, to obtain \eqref{conti1},
we only need to prove that,
for any $f\in C_0(E)$,
\begin{equation}\label{5.3}
  \lim_{\delta\to0}\lim_{n\to\infty} \sup_{t\in[n\delta,(n+1)\delta]}
  \left|\gamma_t\la \phi_0T^{\phi_0}_{(n+1)\delta-t}f,X_t\ra-\la f\phi_0,\psi_0\ra_m W\right|=0.
\end{equation}
Since $\phi_0T^{\phi_0}_{(n+1)\delta-t}f=
e^{-\lambda_0((n+1)\delta-t)}T_{(n+1)\delta-t}(\phi_0f)$,
we only need to show that,
for any $f\in C_0(E)$,
\begin{equation}\label{5.2}
  \lim_{\delta\to0}\lim_{n\to\infty} \sup_{t\in[n\delta,(n+1)\delta]}
  \left|\gamma_t\la T_{(n+1)\delta-t}f,X_t\ra-\la f,\psi_0\ra_m W\right|=0.
\end{equation}
By \eqref{martep2}, we have that, for  any $t\in[n\delta,(n+1)\delta]$,
\begin{align*}
 &\la T_{(n+1)\delta-t}f,X_t\ra=\\
 &=\la T_{\delta}f, X_{n\delta} \ra+\int_{n\delta}^t
 \int_{\mathcal{D}_{<1}(s)}
 \la T_{(n+1)\delta-s}f,\nu\ra \widetilde{N}({\rm d}s,{\rm d}\nu)+\\
 &\qquad+\int_{n\delta}^t
 \int_{\mathcal{D}_{\ge 1}(s)}
 \la T_{(n+1)\delta-s}f,\nu\ra
 \widetilde{N}({\rm d}s,{\rm d}\nu)+
 \int_{n\delta}^t\int_E T_{(n+1)\delta-s}f(x)S^C({\rm d}s,{\rm d}x)\\
 &=:\la T_{\delta}f, X_{n\delta} \ra+H_{t}^{n,\delta}(f)+L_{t}^{n,\delta}(f)+C_{t}^{n,\delta}(f).
\end{align*}
It follows from \eqref{T_t} that
$\phi_0(x)^{-1}T_{\delta}f(x)\in \mathcal{B}_b^+(E).$
Thus, by Theorem \ref{Thm:main1}, we have
$$
\lim_{n\to\infty} \gamma_{n\delta}\la T_{\delta}f, X_{n\delta} \ra
=\la T_{\delta}f,\psi_0\ra_m W
=e^{\lambda_0\delta}\la f,\psi_0\ra_m W.
$$
Note that $\gamma_{(n+1)\delta}\le\gamma_t\le \gamma_{n\delta}$.
Thus,
$$\lim_{\delta\to0}\lim_{n\to\infty} \sup_{t\in[n\delta,(n+1)\delta]}
\left|\gamma_t\la T_{\delta}f, X_{n\delta} \ra-\la f,\psi_0\ra_m W\right|=0.$$

To finish the proof, it suffices to show that
\begin{eqnarray}
  \lim_{n\to\infty} \sup_{t\in[n\delta,(n+1)\delta]}\gamma_{(n+1)\delta}|H_{t}^{n,\delta}(f)|=0;\label{4.12}\\
  \lim_{n\to\infty} \sup_{t\in[n\delta,(n+1)\delta]}\gamma_{(n+1)\delta}|L_{t}^{n,\delta}(f)|=0;\label{4.13}
\end{eqnarray}
and
\begin{equation}\label{4.14}
  \lim_{n\to\infty} \sup_{t\in[n\delta,(n+1)\delta]}\gamma_{(n+1)\delta}|C_{t}^{n,\delta}(f)|=0.
\end{equation}
Using the Markov property of the superprocess $X$, we get that
\begin{align*}
 &\gamma^2_{(n+1)\delta}\P_{\mu}\left(\sup_{t\in[n\delta,(n+1)\delta]}(H_{t}^{n,\delta}(f))^2|\mathcal{G}_{(n-1)\delta}\right)\\
  &= \gamma^2_{(n+1)\delta}\P_{X_{(n-1)\delta}}\left(\sup_{t\in[\delta,2\delta]}
  \left[\int_{\delta}^t
  \int_{\mathcal{D}_{<1}(s+(n-1)\delta)}
  \la T_{2\delta-s}f,\nu\ra \widetilde{N}({\rm d}s,{\rm d}\nu)\right]^2\right)\\
  &\le 4\gamma^2_{(n+1)\delta}\P_{X_{(n-1)\delta}}\left(\left[\int_{\delta}^{2\delta}
  \int_{\mathcal{D}_{<1}(s+(n-1)\delta)}
  \la T_{2\delta-s}f,\nu\ra \widetilde{N}({\rm d}s,{\rm d}\nu)\right]^2\right)\\
  &=4\P_{\mu}\left([H_{(n+1)\delta}(f)-\P_{\mu}(H_{(n+1)\delta}(f)|\mathcal{G}_{n\delta})]^2|\mathcal{G}_{(n-1)\delta}\right),
\end{align*}
where the second to the last line follows from the fact that
$$
\left(\int_{\delta}^t\int_{\mathcal{D}_{<1}(s+(n-1)\delta)}
\la T_{2\delta-s}f,\nu\ra \widetilde{N}({\rm d}s,{\rm d}\nu), t\in[\delta,2\delta]\right)
$$
is a martingale.
Therefore, by \eqref{sqr:H}, we have
$$
\sum_{n=1}^\infty\gamma_{(n+1)\delta}^2
\P_{\mu}\left(\sup_{t\in[n\delta,(n+1)\delta]}(H_{t}^{n,\delta}(f))^2|\mathcal{G}_{(n-1)\delta}\right)
<\infty.
$$
Using the conditional Borel-Cantelli lemma, \eqref{4.12} follows immediately.

Similarly, we can prove that \eqref{4.14} holds. We omit the details here.

Note that
$$|L_{t}^{n,\delta}(f)|\le \int_{n\delta}^{(n+1)\delta}
\int_{\mathcal{D}_{\ge 1}(s)}\la T_{(n+1)\delta-s}f,\nu\ra
\left(N({\rm d}s,{\rm d}\nu)+\widehat{N}({\rm d}s,{\rm d}\nu)\right).$$
Now using \eqref{4.5} and \eqref{4.10} with $m=1$, we immediately get \eqref{4.13}.
The proof is now complete.
\hfill $\Box$

\begin{theorem}\label{LLN}
Suppose that Assumptions  \ref{assum1}--\ref{assum5} hold.
There exists $\Omega_0\subset\Omega$ of probability one (that is,
$\mathbb P_\mu (\Omega_0)=1$ for every
$\mu\in \mathcal M_F(E)$) such that,
for every $\omega\in\Omega_0$ and for every bounded Borel
function $f$ on $E$ satisfying (a) $|f|\le c\phi_0$ for some $c>0$ and (b)
the set of  discontinuous points of $f$ has zero $m$-measure, we have
$$
\lim_{t\to\infty}\gamma_t\la f,X_t\ra(\omega)=\la f,\psi_0\ra_m W(\omega).
$$
\end{theorem}
\noindent\textbf{Proof:}
With the preparation above, the proof of this theorem  is similar to that of \cite[Theorem 1.4]{CRSZ}.
We omit the details here.
\hfill $\Box$

\subsection{Continuous times: Case II}

In this subsection, we will consider the almost sure limit of
$\gamma_t\langle f,  X_t\rangle$
with $f$ being a  general bounded continuous function for some class of superdiffusions.
The underlying process $\xi$ is a diffusion satisfying the following conditions.

Suppose that $E$ is a domain of finite Lebesgue measure in $\R^d$.
Denote by $C_b^1(E)$ the family of  bounded differentiable functions on $E$ whose first order partial derivatives are all continuous.
The underling process $\{\xi,\Pi_x\}$ is a killed diffusion process on $E$ corresponding to the infinitesimal generator
\begin{align}\label{def:L}
L=\frac{1}{2}\nabla\cdot a\nabla+b\cdot \nabla,
\end{align}
where $a$ and $b$ satisfy the following conditions:
\begin{enumerate}
  \item[(a)] $a_{ij}\in C_b^1(E)$, $i,j=1,2\cdots,d,$ and that the matrix $a=(a_{ij})$ is symmetric
  which satisfying, for all $x\in E$ and $v\in \R^d$,
  $$c_0|v|^2\le \sum_{i,j}a_{ij}v_iv_j,$$
  for some positive constant $c_0$.
  \item [(b)]$b_j\in \mathcal{B}_b(E)$, $j=1,\cdots,d$.
\end{enumerate}

Using an argument similar to that in \cite[section 3.2]{CZ}, one can easily show
that $P_t$ has a bounded continuous and strictly positive density $p(t,x,y)$.
Thus Assumption \ref{assum1} holds immediately.
Since $m(E)<\infty$ and the first eigenfunction $\widetilde{\phi}_0\in L^2(E,m)$, we have that $\widetilde{\phi}_0\in L^1(E,m)$.
Then using the fact that $p(1,x,y)$ is bounded and $\widetilde{\phi}_0(x)
=e^{-\widetilde\lambda_0} P_1\widetilde{\phi}_0(x)$, we get that
 $\widetilde{\phi}_0$ is bounded on $E$. Similarly, $\widetilde{\psi}_0$ is also bounded,
which shows that Assumption \ref{assum3}(i) holds. We assume that the semigroup $P_t$ is intrinsically ultracontractive.

 Let $f\in \mathcal{B}_b(E)$, and
 $U^qf$, $q>0$, be the $q$-potential of $f$,
 that is
 $$U^qf(x)=\int_0^\infty e^{-qs}T^{\phi_0}_sf(x){\rm d}s.$$
 For any $t>0$,
 \begin{equation}\label{Uq}
   e^{-qt}T_t^{\phi_0}(U^qf)(x)=\int_t^\infty e^{-qs}T^{\phi_0}_sf(x){\rm d}s.
 \end{equation}

\begin{theorem}\label{thrm:potential}
For any $q>0$, $\mu\in \mathcal{M}_F(E)$ and $f\in \mathcal{B}_b^+(E)$, we have
\begin{equation}\label{potential}
\lim_{t\to\infty}\gamma_t\la \phi_0qU^qf,X_t\ra=\la f\phi_0,\psi_0\ra_m W,\qquad \mbox{a.s.-}\P_{\mu}.
\end{equation}
\end{theorem}
\noindent\textbf{Proof:}
First, we claim that
\begin{equation}\label{claim2}
\lim_{\delta\to0}\lim_{n\to\infty} \sup_{t\in[n\delta,(n+1)\delta]}
\left|\gamma_t\la \phi_0T^{\phi_0}_{(n+1)\delta-t}(U^qf),X_t\ra- \gamma_t\la \phi_0U^qf,X_t\ra\right|=0.
\end{equation}
For any $r\in [0,\delta]$, we have
\begin{align*}
  &q\left|T^{\phi_0}_{r}(U^qf)(x)-U^qf(x)\right|\\
  &=q\left|(e^{qr}-1)\int_{r}^\infty e^{-qs}T^{\phi_0}_sf(x){\rm d}s-\int_0^r e^{-qs}T^{\phi_0}_sf(x){\rm d}s\right|\\
  &\le \|f\|_\infty \left(q(e^{qr}-1)\int_{r}^\infty e^{-qs}{\rm d}s+q\int_0^r e^{-qs}{\rm d}s\right)\\
  &=2\|f\|_\infty (1-e^{-qr}).
\end{align*}
Thus,
\begin{align*}
   &\sup_{t\in[n\delta,(n+1)\delta]}
   \left|\gamma_t\la \phi_0T^{\phi_0}_{(n+1)\delta-t}(qU^qf),X_t\ra- \gamma_t\la \phi_0qU^qf,X_t\ra\right|\\
  &\le   2\|f\|_\infty(1-e^{-q\delta})\sup_{t\in[n\delta,(n+1)\delta]}\gamma_t\la\phi_0,X_t\ra.
\end{align*}
Letting $n\to\infty$ and then $\delta\to0$, the claim \eqref{claim2} follows immediately.
Note that
$$\la qU^qf\phi_0,\psi_0 \ra_m=\la f\phi_0,\psi_0\ra_m.$$
Thus, applying \eqref{5.3} with $f$ replaced by $qU^qf$, we get
\begin{equation}\label{5.4}
  \lim_{\delta\to0}\lim_{n\to\infty} \sup_{t\in[n\delta,(n+1)\delta]}\left|\gamma_t\la \phi_0T^{\phi_0}_{(n+1)\delta-t}(qU^qf),X_t\ra-
  \la \phi_0f,\psi_0\ra_m W\right|=0.
\end{equation}
Now combining \eqref{claim2} and \eqref{5.4}, \eqref{potential} follows immediately.
\hfill $\Box$

\bigskip

\begin{theorem}\label{LLN2}
Suppose $X$ is a superdiffusion on a domain $E\subset\R^d$ of finite Lebesgue measure with spatial motion being a killed diffusion in $E$ with
generator $L$ given in (\ref{def:L}) satisfying conditions (a) and (b).
Suppose that Assumptions \ref{assum1}-\ref{assum5a} hold.
Then there exists $\Omega_0\subset\Omega$ of probability one (that is,
$\mathbb P_\mu (\Omega_0)=1$ for every
$\mu\in \mathcal M_F(E)$) such that,
for every $\omega\in\Omega_0$ and for every bounded Borel
function $f$ on $E$ satisfying (a) $|f|\le c\phi_0$ for some $c>0$ and (b)
the set of  discontinuous points of $f$ has zero $m$-measure, we have
$$
\lim_{t\to\infty}\gamma_t\la f,X_t\ra(\omega)=\la f,\psi_0\ra_m W(\omega).
$$
\end{theorem}
\noindent\textbf{Proof:}
With the preparation above,
 the proof of this theorem is similar to that of \cite[Theorem 1.1]{LRS2}.
We omit the details here.
\hfill $\Box$

\section{Concluding remarks}
Suppose that $X=\{X_t, t\ge 0; \mathbb{P}_{\mu}\}$ is a supercritical superprocess in
a locally compact separable metric space $E$ such that the generator of the mean semigroup
of $X$ has discrete spectrum.
Let $\phi_0$ be a positive
 eigenfunction corresponding to the first eigenvalue $\lambda_0$ of
the generator of the mean semigroup of $X$. Then
$M_t:=e^{-\lambda_0t}\langle\phi_0, X_t\rangle$ is a positive martingale.
Let $M_\infty$ be the limit of $M_t$. It is known
(see \cite{LRS09}) that
$M_\infty$ is non-degenerate iff  the $L\log L$ condition is  satisfied.
In this paper, we prove that, under some further conditions,  there exist
a positive function $\gamma_t$ on $[0, \infty)$
and  a non-degenerate random variable
$W$ such that for any
finite nonzero Borel measure $\mu$ on $E$,
\begin{equation*}
\lim_{t\to\infty}\gamma_t\langle \phi_0,X_t\rangle
=W,\qquad\mbox{a.s.-}\mathbb{P}_{\mu}.
\end{equation*}
We also give the almost sure limit of $\gamma_t\langle f,X_t\rangle$
for a class of general test functions $f$.

In \cite{RSZ19}, a sequel to the present paper, we studied properties of the limit random variable $W$, such as absolute continuity and tail probabilities.

It would be interesting to extend the results of this paper and \cite{RSZ19} to supercritical superprocesses with immigration.

The assumptions of this paper, particularly Assumption \ref{assum3}.(ii), are pretty strong.
For example, supercrtical super Brownian motion in $\R^d$ does not satisfy Assumption \ref{assum3}.(ii).
It would be interesting to consider  corresponding results of
this paper and \cite{RSZ19} for supercritical superprocesses under weaker conditions. It would be very interesting to get rid of the assumption that the generator of the mean semigroup
of $X$ has discrete spectrum.

\bigskip
\noindent
{\bf Acknowledgements:} We thank the referee for very careful reading of the paper and for very helpful comments.

\begin{singlespace}

\end{singlespace}

\vskip 0.2truein
\vskip 0.2truein

\noindent{\bf Yan-Xia Ren:} LMAM School of Mathematical Sciences \& Center for
Statistical Science, Peking
University,  Beijing, 100871, P.R. China. Email: {\texttt
yxren@math.pku.edu.cn}

\smallskip
\noindent {\bf Renming Song:} Department of Mathematics,
University of Illinois,
Urbana, IL 61801, U.S.A
and School of Mathematical Sciences, Nankai University, Tianjin 300071, P. R. China.
Email: {\texttt rsong@illinois.edu}

\smallskip

\noindent{\bf Rui Zhang:} School of Mathematical Sciences, Capital Normal
University,  Beijing, 100048, P.R. China. Email: {\texttt
zhangrui27@cnu.edu.cn}

\end{document}